\documentclass[aos,MSNbibl,nameyear,rotating,dvips]{arximspdf}

\doi{10.1214/14-AOS1214} 
\volume{42}
\issue{3}
\pubyear{2014}
\firstpage{1070}
\lastpage{1101}

\makeatletter
\newtheorem{theorem}{Theorem}[section]
\newtheorem{lemma}{Lemma}[section]
\newtheorem{proposition}{Proposition}[section]
\newproclaim{remark}{Remark}[section]
\newproclaim{ass}{Assumption}
\makeatother

\begin{document}
\begin{frontmatter}

\title{Further results on controlling the false discovery~proportion}
\runtitle{Further results on controlling the FDP}

\begin{aug}
\author[a]{\fnms{Wenge} \snm{Guo}\corref{}\thanksref{t1}\ead[label=e1]{wenge.guo@njit.edu}},
\author[b]{\fnms{Li} \snm{He}\thanksref{t2}\ead[label=e2]{heli@temple.edu}}
\and
\author[b]{\fnms{Sanat K.} \snm{Sarkar}\thanksref{t3}\ead[label=e3]{sanat@temple.edu}}
\runauthor{W. Guo, L. He and S. K. Sarkar}

\affiliation{New Jersey Institute of Technology, Temple University and
Temple University}
\address[a]{W. Guo\\
Department of Mathematical Sciences\\
New Jersey Institute of Technology\\
Newark, New Jersey 07102\\
USA\\
\printead{e1}}
\address[b]{L. He\\
S. K. Sarkar\\
Department of Statistics\\
Temple University\\
Philadelphia, Pennsylvania 19122\\
USA\\
\printead{e2}\\
\phantom{E-mail:\ }\printead*{e3}}
\end{aug}
\thankstext{t1}{Supported by NSF Grants DMS-10-06021 and DMS-13-09162.}
\thankstext{t2}{Supported by Merck Research Fellowship.}
\thankstext{t3}{Supported by NSF Grants DMS-10-06344 and DMS-13-09273.}

\received{\smonth{9} \syear{2011}}
\revised{\smonth{2} \syear{2014}}

%
\begin{abstract}
The probability of false discovery proportion (FDP) exceeding
$\gamma\in[0,1)$, defined as $\gamma$-FDP, has received much
attention as a measure of false discoveries in multiple testing.
Although this measure has received acceptance due to its relevance
under dependency, not much progress has been made yet advancing its
theory under such dependency in a nonasymptotic setting, which
motivates our research in this article. We provide a larger class of
procedures containing the stepup analog of, and hence more powerful
than, the stepdown procedure in Lehmann and Romano [\textit{Ann. Statist.}
\textbf{33} (2005) 1138--1154]
controlling the $\gamma$-FDP under similar positive dependence
condition assumed in that paper. We offer better alternatives of the
stepdown and stepup procedures in Romano and Shaikh
[\textit{IMS Lecture Notes Monogr. Ser.} \textbf{49} (2006a) 33--50,
\textit{Ann. Statist.} \textbf{34} (2006b) 1850--1873]
using
pairwise joint distributions of the null $p$-values. We generalize
the notion of $\gamma$-FDP making it appropriate in situations where
one is willing to tolerate a few false rejections or, due to high
dependency, some false rejections are inevitable, and provide
methods that control this generalized $\gamma$-FDP in two different
scenarios: (i)~only the marginal $p$-values are available and (ii)
the marginal $p$-values as well as the common pairwise joint
distributions of the null $p$-values are available, and assuming
both positive dependence and arbitrary dependence conditions on the
$p$-values in each scenario. Our theoretical findings are being
supported through numerical studies.
\end{abstract}

%

\end{frontmatter}

\section{Introduction}\label{sec1}
The idea of improving the traditional and often too conservative
notion of familywise error rate (FWER) has been one of the main
motivations behind much of the methodological developments taking
place in modern multiple testing. One particular direction in which
this idea has flourished is generalizing the FWER from its original
definition of the probability of at least one false discovery or a
nonzero fraction of false discoveries to one that allows more, yet
tolerable, number or fraction of false discoveries and developing
procedures that control these generalized error rates. The rationale
behind taking this direction is that in many situations where a
large number of hypotheses are tested one is often willing to
tolerate more than one false discovery, controlling of course too
many of them. Moreover, due to high positive dependency among a
group or groups of \mbox{$p$-}values corresponding to true null hypotheses,
as in microarray experiments where the genes involved in the same
biological process or pathway are highly dependent on each other and
exhibit similar expression patterns, it is extremely unlikely that
exactly one null $p$-value will be significant given that at least
one of them will be significant. In such cases, a procedure
controlling the probability of at least $k$ false discoveries, the
$k$-FWER, for some fixed $k > 1$, or the probability of the false
discovery proportion (FDP) exceeding $\gamma$, the $\gamma$-FDP, for
some fixed $0 < \gamma<1$, will have a better ability to detect
more false null hypotheses than the corresponding FWER procedure
(i.e., when $k = 1$ or $\gamma=0$).

Thus, the consideration of the $k$-FWER or $\gamma$-FDP seems more
relevant than that of the FWER when controlling false discoveries in
multiple testing of a large number of hypotheses under dependency.
In fact, it has been noted that the dependency gets naturally
factored into the constructions of procedures controlling the
$k$-FWER or $\gamma$-FDP. For instance, the $k$-dimensional joint
distributions of the null $p$-values can be explicitly used while
constructing procedures controlling the $k$-FWER [\citeauthor{r21} (\citeyear{r21,r22})].
Also, since the FDP becomes more variable and gets more skewed with
increasing dependence among the $p$-values [\citet{r4}, \citet{r11},
\citet{r13}, \citet{r15}, and \citet{r25}], by controlling the tail end
probabilities of the FDP, the $\gamma$-FDP, one considers
controlling a quantity that is more relevant under dependency than
the expected FDP, the false discovery rate (FDR) [\citet{r1}], which is even less conservative than the FWER.

A number of papers have been written over the years on $k$-FWER and
$\gamma$-FDP [\citet{r3}, \citet{r6},
\citet{r7}, \citet{r8}, \citet{r9}, \citet{r12}, \citet{r13}, \citet{r14},
\citeauthor{r16} (\citeyear{r16,r17}),
\citet{r18}, \citet{r19}, \citeauthor{r21} (\citeyear{r21,r22}) and \citet{r27}]. Among these,
\citet{r14}, and \citeauthor{r16} (\citeyear{r16,r17}) are
worth mentioning as they have made some fundamental contributions to
the development of theory and methodology of $\gamma$-FDP. A part of
our research is motivated by these papers, and aims at extending,
and often improving, some results in those papers under certain
dependence situations. The motivation of the other part of our
research comes from the realization that if one indeed is willing to
tolerate a few false rejections, the premise under which one would
seek to use a generalized error rate, the notion of $\gamma$-FDP
does not completely take that into account unless it is further
generalized accordingly. In other words, one should consider in this
case a generalized form of the FDP that accounts for $k$ or more
false rejections, and control the probability of this generalized
FDP, rather than the original FDP, exceeding~$\gamma$. So, we
introduce such a generalized notion of \mbox{$\gamma$-FDP}, called the \mbox{$\gamma
$-}kFDP, and propose
procedures that control it under different dependence scenarios in
this paper.

The paper is organized as follows. We provide some preliminaries in
Section~\ref{sec2}, including the definition of our proposed notion of
$\gamma$-kFDP. Section~\ref{sec3} contains our main results on controlling
the $\gamma$-FDP and $\gamma$-kFDP, developed assuming both positive
dependence (Section~\ref{sec3.1}) and arbitrary dependence (Section~\ref{sec3.2})
conditions on the $p$-values in each of the following two scenarios:
(i) only the marginal $p$-values are available and (ii) the marginal
$p$-values as well as the common pairwise joint distributions of the
null $p$-values are available. We obtain a number of newer results
on $\gamma$-FDP than what are available in the literature. We
construct a larger class of procedures controlling the $\gamma$-FDP
under positive dependence than the stepdown procedure given in
\citet{r14}. This class includes the stepup analog of,
and hence more powerful than, this Lehmann--Romano stepdown
procedure. We offer better alternatives of the stepdown and stepup
procedures in \citeauthor{r16} (\citeyear{r16,r17}), given pairwise joint
distributions of the null $p$-values. Most of our main results have
been obtained through a general framework that allows us not only to
develop procedures controlling the newly proposed notion of
$\gamma$-kFDP, for $k \ge1$, but also to produce the aforementioned
new results on $\gamma$-FDP by taking $k=1$. The performances of the
proposed $\gamma$-FDP and $\gamma$-kFDP procedures, individual as
well as relative to relevant competitors, are numerically
investigated through extensive simulations and reported in Section~\ref{sec4}.
Concluding remarks are made in Section~\ref{sec5}. Proofs of some
supporting results are given in the \hyperref[sec6]{Appendix}.

The supplementary material [\citet{r71}] is added due to space constraints
to include some additional figures related to the numerical
investigations in Section~\ref{sec4}.
Also presented in this section are the findings of simulation studies
conducted to examine the effect of $k$ on a $\gamma$-kFDP controlling
procedure (see Remark~\ref{re2.1}) and to provide an insight into the choice of
$k$ under varying dependence.

\section{Preliminaries}\label{sec2}
Suppose that $H_i, i = 1, \ldots, n$, are the $n$ null hypotheses to
be tested based on their respective $p$-values $P_i, i = 1, \ldots,
n$. Let $P_{(1)} \le\cdots\le P_{(n)}$\vadjust{\goodbreak} be the ordered versions of
all the $p$-values and $H_{(1)}, \ldots, H_{(n)}$ be their
corresponding null hypotheses. There are $n_0$ null hypotheses that
are true. For notational convenience, the $p$-values corresponding
to these true null hypotheses will be denoted by $\widehat P_i, i = 1,
\ldots, n_0$, and their ordered versions by $\widehat P_{(1)} \le\cdots
\le\widehat P_{(n_0)}$.

Multiple testing is typically carried out using a stepwise or
single-step procedure. Given a nondecreasing set of critical values
$0 < \alpha_1 \le\cdots\le\alpha_n < 1$, a stepdown procedure
rejects the set of null hypotheses $\{H_{(i)}, i \le i^*_{\mathrm{SD}}\}$,
where $i^*_{\mathrm{SD}} = \max\{ 1 \le i \le n\dvtx  P_{(j)} \le\alpha_j\
\forall j \le i \}$ if the maximum exists, otherwise accepts all the
null hypotheses. A stepup procedure, on the other hand, rejects the
set of null hypotheses $\{H_{(i)}, i \le i^*_{\mathrm{SU}}\}$, where
$i^*_{\mathrm{SU}} = \max\{ 1 \le i \le n\dvtx  P_{(i)} \le\alpha_i \}$ if the
maximum exists, otherwise accepts all the null hypotheses. A
stepdown or stepup procedure with the same critical values is referred to
as a single-step procedure.

Let $V$ be the number of falsely rejected and $R$ be the total
number of rejected null hypotheses. Then, with $V/R$, which is zero
if $R=0$, defining the false discovery proportion (FDP), and given a
fixed $\gamma\in(0,1)$, the $\gamma$-FDP is defined as the
probability of the FDP exceeding $\gamma$; that is, $\gamma\mbox{-FDP} =
\operatorname{Pr}  (\mathrm{FDP} > \gamma )$. Its generalized version
introduced in this paper, which we call $\gamma$-kFDP, is defined as
follows: let
\[
\mathrm{kFDP} = \cases{ \displaystyle \frac{V}{R}, &\quad if $V \ge k$,
\vspace*{5pt}\cr
0, &\quad otherwise.}
\]
Then $\gamma\mbox{-kFDP} =\operatorname{Pr} (\mbox{kFDP} > \gamma)$. Since
$\gamma$-kFDP is $0$, and hence trivially controlled, for any
procedure if $n_0 < k$, we assume throughout the paper that $k \le n_0
\le n$ when controlling this error rate. Also, while constructing a
$\gamma$-kFDP controlling stepwise procedure, we will consider the
first $k-1$ critical constants to be the same as the $k$th one, as in
$k$-FWER procedures, since their choice does not matter in calculating
the $\gamma$-kFDP.

\begin{remark}\label{re2.1}
It should be noted that since $V$ and FDP are likely
to be highly correlated the distribution of kFDP may be very similar to
that of FDP with a small portion of its lower tail set to $0$.
Therefore, the difference between $\gamma$-kFDP and $\gamma$-FDP may be
realized, with the control over $\gamma$-kFDP providing the stipulated
power improvement, only when $k/n$ exceeds a certain value. Of course,
this value, given a specified $\gamma$, would depend on the type and
strength of dependence. We did a numerical study to verify this
intuition and offer an insight into the choice of $k$ under different
types and varying strengths of dependence, and report its findings in
the supplementary material [\citet{r71}].
\end{remark}

The following is the basic assumption regarding
the marginal distributions of the $p$-values made throughout the
paper.

\begin{ass}\label{ass1}
$\widehat P_i \sim U(0, 1)$.\vadjust{\goodbreak}
\end{ass}

\section{Main results}\label{sec3}
In this section, we present the developments of our stepwise
procedures controlling the $\gamma$-FDP and the newly proposed
$\gamma$-kFDP under both positive dependence and arbitrary
dependence conditions on the $p$-values. Typically, only the
marginal distributions of the null $p$-values are used when
constructing multiple testing procedures. However, in practice, the
null $p$-values often have a known common pairwise joint
distribution, and it would be worthwhile to consider developing
$\gamma$-FDP or $\gamma$-kFDP stepwise procedures explicitly
utilizing such additional dependence information, which could
potentially produce more powerful procedures than just using the
marginal $p$-values. With that in mind, we construct our procedures
in the following two different scenarios under each dependence
condition: (i) only the marginal $p$-values are available, and (ii)
the marginal $p$-values as well as the common pairwise joint
distributions of the null $p$-values are available.

\subsection{Procedures under positive dependence}\label{sec3.1} We will make one of
the following two commonly used assumptions characterizing a
positive dependence structure among the $p$-values.

{\renewcommand{\theass}{2(\textup{a})}
\begin{ass}\label{ass2a}
The conditional
expectation $E  \{\phi( P_1, \ldots, P_{n}) |\break  \widehat P_i \le u
\}$ is nondecreasing in $u \in(0, 1)$ for each $\widehat P_i$
and any nondecreasing (coordinatewise) function $\phi$.
\end{ass}}%

{\renewcommand{\theass}{2(\textup{b})}
\begin{ass}\label{ass2b}
The conditional expectation
$E  \{\phi( \widehat P_1, \ldots, \widehat P_{n_0}) |\break  \widehat P_i \le u
\}$ is nondecreasing in $u \in(0, 1)$ for each $\widehat P_i$
and any nondecreasing (coordinatewise) function $\phi$.
\end{ass}}%

Assumption~\ref{ass2a} is\vspace*{1pt} slightly weaker than that characterized by the
property: $E  \{\phi(P_1, \ldots, P_n) | \widehat P_i = u
\} \uparrow u \in(0, 1)$, referred to as the positive regression
dependence on subset (PRDS) (of the null $p$-values); see, for
example, \citet{BY01}
or \citet{Sa02}. Assumption~\ref{ass2b}, less\vspace*{1pt} restrictive than Assumption~\ref{ass2a}, is a weaker version of the
property: $E  \{\phi( \widehat P_1, \ldots, \widehat P_{n_0}) | \widehat P_i = u  \} \uparrow u \in(0, 1)$, known as the positive
dependence\break  (among the null $p$-values) through stochastic ordering
(PDS) due to \citet{r2}; see also \citet{r23}.

\subsubsection{Based on marginal $p$-values}\label{sec3.1.1}
Under a positive
dependence assumption, \citet{r14} gave a stepdown
procedure controlling the $\gamma$-FDP. We improve this work in two
different ways. First, we consider the stepup analog of this
stepdown procedure, which is known to be always more powerful in the
sense of discovering more, and prove that it also controls the
$\gamma$-FDP under the same assumption. Second, we offer larger
class of stepdown and stepup procedures controlling the $\gamma$-FDP
under similar assumption. The procedures in this larger class are
presented in a general framework allowing us to propose procedures
controlling not only the $\gamma$-FDP but also the $\gamma$-kFDP for
$k \ge2$.

\begin{theorem}\label{th3.1}
The stepup or stepdown procedure with the critical constants
%
\begin{equation}\label{equ1}
\alpha_i = \frac{(\lfloor\gamma i \rfloor+ 1) \alpha}{ n +
\lfloor\gamma i \rfloor+ 1 - i},\qquad i=1, \ldots, n,
\end{equation}
controls the $\gamma$-FDP at $\alpha$ under Assumptions~\ref{ass1} and~\ref{ass2b}.
\end{theorem}

\begin{pf}
Let $g(R)= \lfloor\gamma R \rfloor+ 1$. Then first note
that
%
\begin{eqnarray}\label{equ2}
\bigl\{ V \ge g(R) \bigr\} &=& \bigcup_{v = 1 }^{n_0}
\bigl\{ \widehat P_{(v)} \le\alpha_{R}, g(R) \le v, V=v \bigr\}
\nonumber
\\
& = & \bigcup_{v = 1}^{n_0} \biggl\{ \widehat P_{(v)} \le \frac{g(R)\alpha}{n-R+g(R)}, g(R) \le v, V=v \biggr\}
\nonumber\\[-8pt]\\[-8pt]
& \subseteq& \bigcup_{v = 1}^{n_0} \biggl\{
\widehat P_{(v)} \le \frac{v\alpha}{n-R+v}, V=v \biggr\}\nonumber
\\
& \subseteq& \bigcup_{v = 1}^{n_0} \biggl\{
\widehat P_{(v)} \le \frac{v\alpha}{n_0}, V=v \biggr\}
\subseteq \bigcup _{v =
1}^{n_0} \biggl\{ \widehat P_{(v)} \le
\frac{v\alpha}{n_0} \biggr\}.\nonumber
\end{eqnarray}
The probability of the event in the right-hand side
of (\ref{equ2}) is known to be less than or equal to $\alpha$ under
Assumptions~\ref{ass1} and~\ref{ass2b} from the so-called Simes' inequality [\citet{r26}, \citet{r20}, \citet{r24}]. Thus, we get the
desired result noting that $\gamma\mbox{-FDP} = \operatorname{Pr}  (V \ge
g(R)  )$.
\end{pf}

\begin{remark}\label{re3.1}
\citet{r14} proposed only the stepdown
procedure considered in Theorem~\ref{th3.1} under the same assumptions.
Thus, Theorem~\ref{th3.1} provides an improvement of the Lehmann--Romano
result, since we now have an alternative procedure under the same
assumptions, the stepup one, which is theoretically known to be
more powerful. Moreover, not only our proof of the $\gamma$-FDP
control is much simpler but also it covers both ours and the
Lehmann--Romano original stepdown procedures. Our simulation studies
indicate that this power improvement can be obvious
when the underlying test statistics are highly
correlated (see Figure~\ref{fig1} and Figures S.1--S.3 in the supplementary
material [\citet{r71}]).
\end{remark}

%
\begin{figure}

\includegraphics{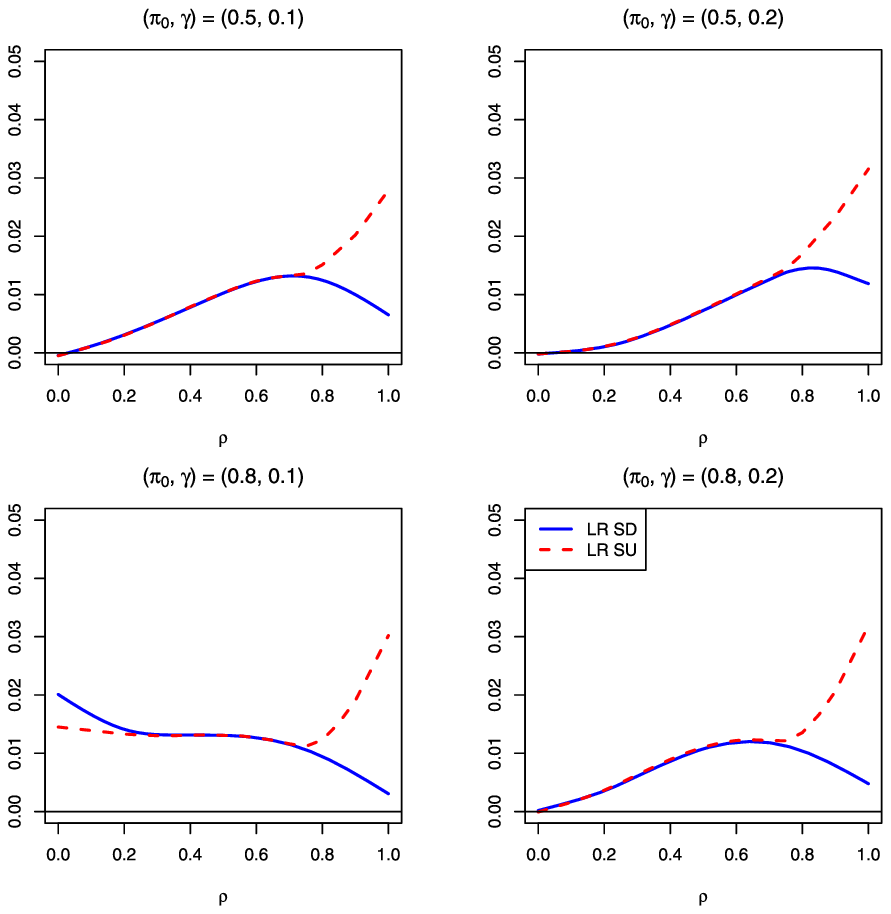}

\scriptsize{(a)~Simulated $\gamma$-FDP}
\caption{Simulated values of $\gamma$-FDP and average power of the original
Lehmann--Romano stepdown procedure (LR SD) and its stepup analogue
(LR~SU), for $n=100$ and $\alpha= 0.05$.}\label{fig1}
\end{figure}

\setcounter{figure}{0}
%
\begin{figure}

\includegraphics{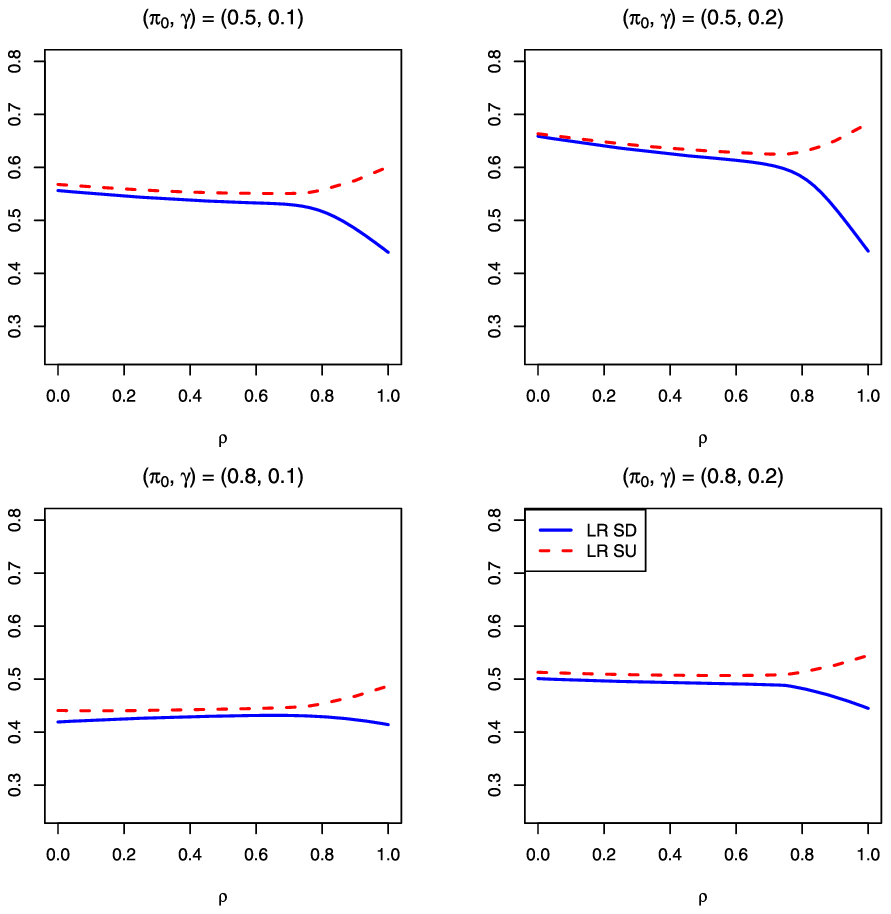}

\scriptsize{(b) Simulated average power}
\caption{(Continued).}
\end{figure}

There are more general results than Theorem~\ref{th3.1} in terms of deriving
procedures controlling the $\gamma$-FDP under Assumptions~\ref{ass1} and~\ref{ass2a}
or~\ref{ass2b}. More specifically, we can start with any stepdown or stepup
procedure, which may or may not control the $\gamma$-FDP to begin
with, and rescale its critical values using a suitable upper bound
for its $\gamma$-FDP derived under Assumptions~\ref{ass1} and~\ref{ass2a} or~\ref{ass2b} so
that the $\gamma$-FDP based on these modified critical values is
ultimately controlled. \citeauthor{r16} (\citeyear{r16,r17}) first developed
this idea, but they did it without any positive dependence
assumption. We are now going to present these results in the general
framework of controlling the $\gamma$-kFDP.

Our next main result is obtained with the idea of constructing a
stepdown procedure controlling the $\gamma$-kFDP under Assumptions~\ref{ass1}
and~\ref{ass2a}. The following lemma, to be proved in the \hyperref[sec6]{Appendix}, will provide
the starting point for the development of this procedure.

\begin{lemma}\label{le3.1} With $n_1=n-n_0$, let $M=\min\{n_0,\lfloor\gamma
n_1/(1-\gamma)
\rfloor+ 1\}$, and $m(i) = \max\{ 0 \le j \le n_1\dvtx  \lfloor\gamma
j/(1-\gamma) \rfloor+ 1 = i \}$, for each $i=1, \ldots, M$, where
$m(0)=0$. Consider a stepdown procedure with critical values $\alpha_1
\le\cdots\le\alpha_n$. Let $S$ be the number of
rejected false null hypotheses. Then
%
\begin{eqnarray}\label{equ3}
&& I \bigl( V > \max[\gamma R, k-1] \bigr)
\nonumber\\[-8pt]\\[-8pt]
&&\qquad \le \sum_{i = 1 }^{M} I \bigl( \widehat P_{(i \vee k)} \le\alpha_{i
\vee k +m(i)}, \bigl\lfloor\gamma S/(1-\gamma) \bigr
\rfloor+1 = i \bigr),\nonumber
\end{eqnarray}
for any fixed $1 \le k \le n_0$.
\end{lemma}

Taking expectations of both sides in (\ref{equ3}), we note that
%
\begin{eqnarray}\label{equ4}
\gamma\mbox{-kFDP} &=& \operatorname{Pr} \bigl\{ V > \max[\gamma R, k-1] \bigr
\}\nonumber
\\
& \le& \sum_{i=1}^{M} \operatorname{Pr} \bigl( \widehat P_{(i \vee k)} \le\alpha _{i \vee k +
m(i)}, \bigl\lfloor\gamma S / (1 - \gamma)
\bigr\rfloor+ 1 = i \bigr)
\\
& \le& \sum_{j = 1}^{n_0} \sum
_{i=1}^{M} \frac{1}{i \vee k}\operatorname{Pr} \bigl( \widehat P_j \le\alpha_{i \vee k + m(i)}, \bigl\lfloor\gamma S / (1 - \gamma)
\bigr\rfloor+ 1 = i \bigr)\nonumber
\\
& = & \sum_{j = 1}^{n_0} \sum
_{i=1}^{M} \frac{\alpha_{i \vee k +m(i)}}{i \vee k} \operatorname{Pr} \bigl( \bigl
\lfloor \gamma S / (1 - \gamma) \bigr\rfloor+ 1 = i | \widehat P_j \le
\alpha_{i \vee k + m(i)} \bigr)
\nonumber
\\
& \le& \max_{1 \le i \le M} \biggl\{ \frac{\alpha_{i \vee k+ m(i)}}{i
\vee k} \biggr\}\nonumber
\\
&&{}\times \sum
_{j = 1}^{n_0} \sum_{i=1}^{M}
\operatorname{Pr} \bigl( \bigl\lfloor \gamma S / (1 - \gamma) \bigr\rfloor+ 1 = i | \widehat P_j \le \alpha_{i \vee k + m(i)} \bigr),
\nonumber
\end{eqnarray}
with the second
inequality following from this:
%
\begin{eqnarray}\label{equ5}
I (\widehat P_{(i)} \le t ) & \le& \frac{1}{i} \sum
_{j=1}^{n_0} I (\widehat P_j \le t )\qquad\mbox{for any constant } 0 < t <1,
\end{eqnarray}
which can be obtained from Markov's inequality.

Now, for each $1 \le j \le n_0$, we have
%
\begin{eqnarray}\label{equ6}
& & \sum_{i=1}^M \operatorname{Pr} \bigl( \bigl
\lfloor\gamma S / (1 - \gamma) \bigr\rfloor+ 1 = i | \widehat P_j \le
\alpha_{i \vee k + m(i)} \bigr)
\nonumber
\\
&&\qquad  = \sum_{i=1}^M\operatorname{Pr} \bigl( \bigl
\lfloor\gamma S / (1 - \gamma) \bigr\rfloor+ 1 \ge i | \widehat P_j \le
\alpha_{i \vee k + m(i)} \bigr)\nonumber
\\
&&\quad\qquad{}- \sum_{i=1}^M\operatorname{Pr} \bigl( \bigl
\lfloor\gamma S / (1 - \gamma) \bigr\rfloor+ 1 \ge i+1 | \widehat P_j \le
\alpha_{i \vee k + m(i)} \bigr)
\nonumber\\[-8pt]\\[-8pt]
&&\qquad  \le \sum_{i=1}^M\operatorname{Pr} \bigl( \bigl
\lfloor\gamma S / (1 - \gamma) \bigr\rfloor+ 1 \ge i | \widehat P_j \le
\alpha_{i \vee k + m(i)} \bigr)
\nonumber
\\
&&\quad\qquad{}- \sum_{i=1}^M\operatorname{Pr} \bigl( \bigl
\lfloor\gamma S / (1 - \gamma) \bigr\rfloor+ 1 \ge i+1 | \widehat P_j \le
\alpha_{(i+1) \vee k + m(i+1)} \bigr)
\nonumber
\\
&&\qquad \le \operatorname{Pr} \bigl( \bigl\lfloor\gamma S / (1 - \gamma) \bigr\rfloor+ 1 \ge1
| \widehat P_j \le\alpha_{k + m(1)} \bigr)
= 1.\nonumber
\end{eqnarray}
The first inequality follows from Assumption~\ref{ass2a}, since $I  (
\lfloor\gamma S / (1 - \gamma) \rfloor+\break  1 \ge i  )$ is a
decreasing function of all the $p$-values. Applying (\ref{equ6}) to (\ref{equ4}), we
finally note
\begin{eqnarray}
\gamma\mbox{-kFDP} & \le& \max_{1 \le i \le M} \biggl\{
\frac{n_0 \alpha_{i \vee k +m(i)}}{i \vee k} \biggr\},
\nonumber
\end{eqnarray}
and thus we have our next main result as follows.

\begin{theorem}\label{th3.2}
Let $M$ and $m(i)$, for $i=1, \ldots, M$, be defined as in Lemma~\ref{le3.1}.
Then,\vspace*{1pt} given any set of constants $0 = \alpha_0^{\prime}
\le\alpha_1^{\prime} \le\cdots\le\alpha_n^{\prime}$, the stepdown
procedure with the critical values $\alpha_{i \vee k} = \alpha
\alpha_{i \vee k}^{\prime}/C_{k,n,\mathrm{SD}}^{(1)}, i=1, \ldots, n$, where
\begin{eqnarray}
C_{k,n,\mathrm{SD}}^{(1)} & = & \max_{k \le n_0 \le n } \max
_{1
\le
i \le M} \biggl\{ \frac{n_0 \alpha_{i \vee k +m(i)}^{\prime}}{i \vee
k} \biggr\},
\nonumber
\end{eqnarray}
controls the $\gamma$-kFDP at $\alpha$ under
Assumptions~\ref{ass1} and~\ref{ass2a}.
\end{theorem}

A stepup analog of Theorem~\ref{th3.2} can be developed starting from the
following lemma, whose proof again is given in the \hyperref[sec6]{Appendix}.

\begin{lemma}\label{le3.2}
Let $\tilde m (i) = \min\{ m^*(i), i +n_1\}$, where
$m^*(i) = \max\{1\le j \le n\dvtx  \lfloor\gamma j \rfloor+ 1 \le i \}$,
for each $i=1,
\ldots, n_0$, and $m^*(0)=0$. Consider a stepup\vadjust{\goodbreak} procedure with critical
values $\alpha_1 \le\cdots\le\alpha_n$. Then, for any fixed $1 \le k
\le n_0$,
%
\begin{eqnarray}\label{equ7}
I \bigl(V > \max[\gamma R,k-1] \bigr)
&\le& \sum_{j = 1}^{n_0}
\sum_{i = k}^{n_0} \frac{I ( \widehat P_j \le\alpha
_{\tilde m(i)}, \widehat R_2 = i  ) }{i }\nonumber
\\
&\le& \sum_{j = 1}^{n_0} \frac{I ( \widehat P_j \le\alpha_{\tilde m(k)}, \widehat R_2 \ge k  ) }{k }
\\
&&{}  +
\sum_{j = 1}^{n_0} \sum
_{i = k+1}^{n_0} \frac{ I  ( \alpha_{\tilde m(i-1)} < \widehat P_j
\le\alpha_{\tilde m(i)}, \widehat R_2 \ge i  ) }{ i },\nonumber
\end{eqnarray}
with the double summation in the right-hand side of the second
inequality being zero if $n_0=k$, where $\widehat R_2$ is the number of
rejections in a stepup procedure
based on the $p$-values $\widehat P_i$, $i=1, \ldots, n_0$, and the
critical values $\alpha_{\tilde m(i)}$, $i=1, \ldots, n_0$.
\end{lemma}

\begin{remark}\label{re3.2}
If we let $n_0 = n$ in the above lemma, we note that
$I(V > \max[\gamma R, k-1]) = I(V \ge k)$ and $\tilde m(i) = i$.
In other words, the above lemma produces inequalities similar to (\ref{equ7})
for $I (\widehat R_2 \ge k)$, with $\widehat R_2$ representing the number of
rejections in a stepup procedure based on the null $p$-values $\widehat P_i$, $i=1, \ldots, n_0$, and critical values $\alpha_{i}$, $i=1, \ldots, n_0$. For instance, from the second inequality in~(\ref{equ7}), we have
%
\begin{eqnarray}\label{equ8}
I (\widehat R_2 \ge k ) & \le& \sum_{j =
1}^{n_0}
\frac{I ( \widehat P_j \le\alpha_{k}  ) }{k } + \sum_{j = 1}^{n_0} \sum
_{i = k+1}^{n_0} \frac{ I  ( \alpha_{i-1} <
\widehat P_j \le\alpha_{i}  ) }{ i },
\end{eqnarray}
which will be of use later. Of course, the first inequality in this
case becomes an equality.
\end{remark}

Taking expectations of both sides of the first inequality in (\ref{equ7}), we get
%
\begin{eqnarray}\label{equ9}
\gamma\mbox{-kFDP} & \le& \sum_{j = 1}^{n_0}
\sum_{i = k}^{n_0} \frac{\alpha_{\tilde m(i)}}{ i } \operatorname{Pr} (
\widehat R_2 = i | \widehat P_j \le\alpha_{\tilde m(i)} )
\nonumber\\[-8pt]\\[-8pt]
& \le& \max_{k \le i \le n_0} \biggl\{ \frac{\alpha_{\tilde m(i)}}{ i } \biggr\} \sum
_{j = 1}^{n_0} \sum_{i = k}^{n_0}
\operatorname{Pr} ( \widehat R_2 = i | \widehat P_j \le
\alpha_{\tilde m(i)} ).\nonumber
\end{eqnarray}
Making the same kind of arguments as in (\ref{equ6}), we note that
\begin{eqnarray}
& & \sum_{i = k}^{n_0} \operatorname{Pr} ( \widehat R_2 = i | \widehat P_j \le\alpha_{\tilde m(i)} ) \le\operatorname{Pr} ( \widehat R_2 \ge k |\widehat P_j \le\alpha_{\tilde m(k)}
) \le1,
\nonumber
\end{eqnarray}
for each $1 \le j \le n_0$, using
the fact that $I  (\widehat R_2 \ge i  )$ is a decreasing
function of the null $p$-values and applying Assumption~\ref{ass2b}. Hence,
\begin{eqnarray}
\gamma\mbox{-kFDP} & \le& \max_{k \le i
\le n_0} \biggl\{
\frac{ n_0 \alpha_{\tilde m(i)}}{ i } \biggr\},
\nonumber
\end{eqnarray}
which provides the following result.

\begin{theorem}\label{th3.3}
Let ${\tilde m}(i)$ be defined as in Lemma~\ref{le3.2} for
$i=1, \ldots, n_0$. Then, given any set\vspace*{1pt} of constants $0
=\alpha_0^{\prime} \le\alpha_1^{\prime} \le\cdots\le
\alpha_n^{\prime}$, the stepup procedure with the critical
values $\alpha_{i \vee k} = \alpha\alpha_{i \vee k}^{\prime
}/C_{k,n,\mathrm{SU}}^{(1)}, i=1,
\ldots, n$, where
\[
C_{k,n,\mathrm{SU}}^{(1)} = \max_{k \le n_0
\le n } \max
_{k \le i \le n_0} \biggl\{\frac{ n_0 \alpha_{\tilde m(i)}^{\prime}}{ i } \biggr\},
\]
controls the $\gamma$-kFDP at $\alpha$ under Assumptions~\ref{ass1} and~\ref{ass2b}.
\end{theorem}

\begin{remark}\label{re3.3}
Theorems~\ref{th3.2} and~\ref{th3.3} not only provide general
approaches to constructing stepdown and stepup $\gamma$-kFDP
controlling procedures, respectively, using only the marginal
$p$-values under independence or certain positive dependence
condition on the $p$-values, but also produce results when $k=1$
that improve some previous works on controlling the $\gamma$-FDP
[\citet{r14}, \citeauthor{r16} (\citeyear{r16,r17})]. For
instance, if we choose the $\alpha_i^{\prime}$'s in these theorems
as follows: $\alpha_{i}^{\prime} = \{\lfloor\gamma i \rfloor+1
\}\alpha/\{n +\lfloor\gamma i \rfloor+1 - i\}$, $i=1, \ldots, n$,
then we get the original Lehmann--Romano procedure and its stepup
analog, since both $C_{1, n,\mathrm{SD}}^{(1)}$ and $C_{1, n,\mathrm{SU}}^{(1)}$ are
equal to $\alpha$ (see\vspace*{1pt} Proposition~\ref{pr6.1} and its proof in
the \hyperref[sec6]{Appendix}). However, there are other stepdown and stepup procedures
controlling the $\gamma$-FDP under these assumptions, such as those
obtained by re-scaling the critical values, $\alpha_i^{\prime} = i
\alpha/n$, $i=1, \ldots, n$, of the BH [\citet{r1}] stepup or the critical values, $\alpha_i^{\prime} = i\alpha/[n
- i(1- \alpha) +1]$, $i=1, \ldots, n$, of the GBS [\citet{r5}] stepdown methods, that can be
constructed using the above theorems. Our simulation studies
indicate that a stepwise procedure based on the rescaled versions
of the BH or GBS critical values is less powerful than that based on
the rescaled version of the Lehmann--Romano critical values (see
Figures~S.4--S.7 in the supplementary material [\citet{r71}]). Therefore, the interest
of Theorems~\ref{th3.2} and~\ref{th3.3} with respect to Theorem~\ref{th3.1} seems to be
mainly theoretical when $k=1$.
\end{remark}

\subsubsection{Based on marginal and pairwise null distributions of the $p$-values}\label{sec3.1.2}
In practice, the null $p$-values often have a known common pairwise joint
distribution, and by explicitly
utilizing such correlation information better adjustments can be
made, potentially resulting in more powerful $\gamma$-FDP stepwise
procedures. So, with that in mind, we present some results here and
in Section~\ref{sec3.2.2} under the following assumption along with
Assumptions~\ref{ass1} and~\ref{ass2b} or only Assumption~\ref{ass1}.

\setcounter{ass}{2}
\begin{ass}\label{ass3}
The null $p$-values
$\widehat P_1, \ldots, \widehat P_{n_0}$ have a known common pairwise joint
distribution function $F(u,v) = \operatorname{Pr}  (\widehat P_i \le u, \widehat P_j \le v  )$.
\end{ass}

We consider generalizing the Lehmann--Romano stepwise procedure in
Theorem~\ref{th3.1}, for any fixed $2 \le k \le n_0$. The $\gamma$-kFDP of
this procedure is given by
%
\begin{eqnarray}\label{equ10}
\gamma\mbox{-kFDP} & = & \operatorname{Pr} \bigl\{ V \ge\max\bigl[g(R),k\bigr]
\bigr\}
\nonumber
\\
& = & \operatorname{Pr} \Biggl( \bigcup_{v = k }^{n_0}
\bigl\{ \widehat P_{(v)} \le\alpha_{R}, g(R) \le v, V=v \bigr\}
\Biggr)
\\
& \le& \operatorname{Pr} \Biggl(\bigcup_{v = k}^{n_0}
\biggl\{ \widehat P_{(v)} \le\frac{v \alpha}{n_0} \biggr\} \Biggr)
=  \operatorname{Pr} (
\widehat R_{n_0} \ge k ),\nonumber
\end{eqnarray}
where $\widehat R_{n_0}$ is the number of rejections in the stepup procedure based
on all the $n_0$ null $p$-values and the critical values $\beta_i=i
\alpha/n_0$, $i=1, \ldots, n_0$. The $\gamma$-kFDP can be bounded
above using the following inequality which holds under Assumptions~\ref{ass1}
and~\ref{ass2b}:
%
\begin{equation}\label{equ11}
\operatorname{Pr} (\widehat R_{n_0} \ge k ) \le\frac{\alpha}{n_0} \sum
_{i=1}^{n_0} \operatorname{Pr} \bigl( \widehat R_{n_0-1}^{(-i)}
\ge k-1|\widehat P_{i} \le\beta_{k} \bigr),
\end{equation}
for any fixed $1 \le k \le n_0$, where $\widehat R_{n_0-1}^{(-i)}$ is the number of rejections in the stepup
procedure based on the $n_0-1$ null $p$-values $\{\widehat P_1, \ldots,
\widehat P_{n_0} \} \setminus\{\widehat P_i\}$ and the critical values
$\beta_i$, $i=2, \ldots, n_0$. This can be proved using arguments
similar to those used above while proving Theorems~\ref{th3.2} or~\ref{th3.3}; see
the \hyperref[sec6]{Appendix}, for a proof.

As seen from (\ref{equ11}), if we rely only on the marginal distributions of
the null $p$-values, we simply get $\gamma$-kFDP $\le\alpha$, and
thus our attempt to generalize the Lehmann--Romano procedure to a
$\gamma$-kFDP controlling procedure under Assumption~\ref{ass2b} does not
work in the sense that it takes us back to the original
Lehmann--Romano procedure, which is trivially known to control the
$\gamma$-kFDP. Hence, we consider utilizing also the pairwise
distributions of the null $p$-values to obtain a nontrivial
generalization of the Lehmann--Romano procedure. More specifically,
we use the following inequality provided by (\ref{equ8}):
%
\begin{eqnarray}\label{equ12}
& & I \bigl(\widehat R_{n_0-1}^{(-i)} \ge k -1 \bigr)
\nonumber\\[-8pt]\\[-8pt]
&&\qquad  \le \sum_{j (\neq i)= 1}^{n_0} \frac{I ( \widehat P_j \le
\beta_{k}  ) }{k-1 } +
\sum_{j (\neq i)= 1}^{n_0} \sum
_{l =
k}^{n_0-1} \frac{ I  ( \beta_{l} < \widehat P_j \le\beta_{l+1}
) }{ l },\nonumber
\end{eqnarray}
and apply it to the right-hand side
of (\ref{equ11}) to get the following upper bound for the $\gamma$-kFDP of
the Lehmann--Romano stepwise procedure under Assumption~\ref{ass2b}:
%
\begin{eqnarray}\label{equ13}
\gamma\mbox{-kFDP} & \le& \frac{\alpha
}{n_0} \sum
_{i=1}^{n_0} \sum_{j (\neq i)= 1}^{n_0}
\Biggl( \frac{\operatorname{Pr}  ( \widehat P_j \le\beta_{k}|\widehat P_{i} \le\beta_k  )
}{k-1 }
\nonumber\\[-8pt]\\[-8pt]
&&\hspace*{67pt}{} +   \sum_{l = k}^{n_0-1}
\frac{ \operatorname{Pr}  ( \beta_{l} < \widehat P_j \le\beta_{l+1}|\widehat P_i
\le\beta_k  ) }{ l } \Biggr).\nonumber
\end{eqnarray}

Based on this upper bound and that the $\gamma$-kFDP of
the Lehmann--Romano stepwise procedure is $\le\alpha$ under
Assumptions~\ref{ass1} and~\ref{ass2b}, we now have the following theorem providing the
desired generalized version of the Lehmann--Romano procedure
controlling the $\gamma$-kFDP.

\begin{theorem}\label{th3.4}
Let $2 \le k \le n_0$ and Assumption~\ref{ass3} hold. Given
$\beta_i = i \alpha/n_0$, $i=1, \ldots, n_0$, let
%
\begin{equation}\label{equ14}
\qquad C_{k,n} = \max_{k \le n_0 \le n} \Biggl\{ (n_0 -1)
\Biggl(\frac{F(\beta
_k|\beta_k)}{k-1} + \sum_{l=k}^{n_0-1}
\frac{F(\beta_{l+1}|\beta_k) -
F(\beta_{l}|\beta_k)}{l} \Biggr) \Biggr\}
\end{equation}
with the summation within parentheses being zero if $n_0=k$, where
$F(u|v) = F(u,v)/v$. Then the stepup or stepdown procedure with the
critical constants $\alpha_{i \vee k}$, $i=1, \ldots, n$, where
%
\begin{equation}\label{equ15}
\alpha_i = \frac{(\lfloor\gamma i \rfloor+
1)\alpha}{(C_{k,n} \wedge1) (n + \lfloor\gamma i \rfloor+ 1 -
i)},\qquad i=1, \ldots, n,
\end{equation}
controls the $\gamma$-kFDP
at $\alpha$ under Assumptions~\ref{ass1} and~\ref{ass2b}.
\end{theorem}

\subsection{Procedures under arbitrary dependence}\label{sec3.2}
We now present some $\gamma$-kFDP controlling procedures under
arbitrary dependence of the $p$-values. By arbitrary dependence, we
mean that these $p$-values are not known to have any specific type
of dependence structure, like positive or other, even though their
joint distributions of some particular orders might be known. We
will assume, as in Section~\ref{sec3.1.2}, that the null $p$-values have a
common pairwise joint distribution of a known form $F(u,v)$. Our
procedures are developed relying either only on the marginal
$p$-values or on the marginal as well as this common pairwise joint
null distribution of the $p$-values. We can obtain some new results
on controlling the $\gamma$-FDP by taking $k=1$.

\subsubsection{Based on marginal $p$-values}\label{sec3.2.1}

First, let us consider a stepdown procedure with critical values
$\alpha_{i \vee k}$, $i=1, \ldots, n$. Starting from Lemma~\ref{le3.1} and
proceeding as in proving Theorem~\ref{th3.2}, we have, with $i \vee k +
m(i)$ defined as ${\bar m}(i)$ [where ${\bar m }(0) =0$] for
notational convenience,
%
\begin{eqnarray}\label{equ16}
& & I \bigl( V > \max[ \gamma R, k-1] \bigr)
\nonumber
\\
&&\qquad \le \sum_{j=1}^{n_0} \sum
_{i=1}^M \frac{I  (\widehat P_j \le\alpha_{\bar
{m}(i)}  )}{i \vee k} I \bigl(\bigl\lfloor
\gamma S/(1-\gamma) \bigr\rfloor+1 = i \bigr)
\nonumber\\[-8pt]\\[-8pt]
&&\qquad \le \sum_{j=1}^{n_0} \sum
_{i = 1}^{M} \biggl[ \frac{I
(\widehat P_j \le\alpha_{\bar{m}(i)}  )}{i \vee k} -
\frac{I
(\widehat P_j
\le\alpha_{\bar{m}(i-1)}  )}{(i-1) \vee k} \biggr] I \bigl(\bigl\lfloor\gamma S/(1-\gamma) \bigr\rfloor+1
\ge i \bigr)\nonumber\hspace*{-15pt}
\\
&&\qquad \le \sum_{j=1}^{n_0} \sum
_{i = 1}^{M} \frac{I
(\alpha_{\bar{m}(i-1)}< \widehat P_j \le\alpha_{\bar{m}(i)}  )}{i
\vee k}.
\nonumber
\end{eqnarray}
Taking expectations of both sides in (\ref{equ16}), we get
%
\begin{eqnarray}\label{equ17}
\gamma\mbox{-kFDP} & = & \operatorname{Pr} \bigl\{ V \ge\max\bigl[ g(R), k\bigr]
\bigr\}
\le n_0 \sum_{i=1}^{M}
\frac{\alpha_{\bar{m}(i)}-
\alpha_{\bar{m}(i-1)}}{ i \vee k}.
\end{eqnarray}
From this, we get the following theorem.

\begin{theorem}\label{th3.5}
Let $M$ and $m(i)$, for $i=1, \ldots, M$, be defined as in Lemma~\ref{le3.1},
and ${\bar m}(i) = i \vee k + m(i)$ [where ${\bar m }(0) =0$].
Then, given any set of constants $\alpha_k^{\prime} \le\cdots\le
\alpha_n^{\prime}$, the stepdown\vspace*{1pt} procedure with the critical values
$\alpha_i = \alpha\alpha_{i \vee k}^{\prime}/C_{k,n,\mathrm{SD}}^{(2)},
i=1, \ldots, n$, where
\begin{eqnarray}
C_{k,n,\mathrm{SD}}^{(2)} & = & \max_{k \le n_0 \le n } \Biggl\{
n_0 \Biggl(\sum_{i=1}^{M}
\frac{\alpha_{{\bar m}(i)}^{\prime}- \alpha_{{\bar
m}(i-1)}^{\prime}}{i \vee k} \Biggr) \Biggr\},
\nonumber
\end{eqnarray}
controls the $\gamma$-kFDP at $\alpha$ under Assumption~\ref{ass1}.
\end{theorem}

We now present the development of a stepup analog of Theorem~\ref{th3.5}. From
Lemma~\ref{le3.2}, we note that for a stepup procedure with critical values
$\alpha_{i \vee k}$, $i=1, \ldots, n$,
%
\begin{eqnarray}\label{equ18}
& & I \bigl( V > \max[ \gamma R, k-1] \bigr)
\nonumber\\[-9pt]\\[-9pt]
&&\qquad \le \sum_{j=1}^{n_0} \frac{I  ( \widehat P_j \le
\alpha_{\tilde m(k)}  )}{k}+
\sum_{j=1}^{n_0} \sum
_{i=k+1}^{n_0} \frac{I  (\alpha_{\tilde m (i-1)}< \widehat P_j
\le\alpha_{\tilde m(i)}  )}{i}.\nonumber
\end{eqnarray}
Taking expectations of both sides in (\ref{equ18}), we get
%
\begin{eqnarray}\label{equ19}
\gamma\mbox{-kFDP} & = & \operatorname{Pr} \bigl\{ V \ge\max[\gamma R, k-1]
\bigr\}
\nonumber\\[-9pt]\\[-9pt]
& \le& n_0 \Biggl(\frac{\alpha
_{\tilde m(k)}}{k}+\sum
_{i=k+1}^{n_0} \frac{\alpha_{\tilde m(i)}- \alpha_{\tilde m(i-1)}}{i} \Biggr),\nonumber
\end{eqnarray}
which
gives the following theorem.

\begin{theorem}\label{th3.6}
Let $\tilde m(i)$ be defined as in Lemma~\ref{le3.2} for $i=1, \ldots,
n_{0}$. Then, given any set of constants $\alpha_k^{\prime} \le
\cdots\le\alpha_n^{\prime}$, the stepup procedure with the
critical values $\alpha_i = \alpha\alpha_{i \vee
k}^{\prime}/C_{k,n,\mathrm{SU}}^{(2)}, i=1, \ldots, n$, where
\[
C_{k,n,\mathrm{SU}}^{(2)} = \max_{k \le n_0 \le n } \Biggl\{
n_0 \Biggl( \frac{\alpha_{\tilde m(k)}^{\prime}}{k} +\sum_{i=k+1}^{n_0}
\frac{\alpha_{\tilde m(i)}^{\prime} -
\alpha_{\tilde m(i-1)}^{\prime}}{i} \Biggr) \Biggr\},
\]
controls the $\gamma$-kFDP at $\alpha$ under
Assumption~\ref{ass1}.
\end{theorem}

%
\begin{remark}\label{le3.4}
When $k=1$, the results in Theorems~\ref{th3.5} and~\ref{th3.6} reduce to those
given by \citeauthor{r16}  in (\citeyear{r16}) and (\citeyear{r17}), respectively,
although\vadjust{\goodbreak} our expressions of the upper bounds given in these theorems
are different from theirs. Thus, our results generalize those of
Romano and Shaikh from controlling the $\gamma$-FDP to $\gamma$-kFDP
under arbitrary dependence and relying only on the marginal null
distributions of the $p$-values. However, we should emphasize that
we provide alternative, much simpler proofs for these results.
\end{remark}

\subsubsection{Based on marginal and pairwise distributions of the null $p$-values}\label{sec3.2.2}

We will start again from Lemma~\ref{le3.1} towards constructing a stepdown
procedure. Consider splitting the sum in the right-hand side of (\ref{equ3})
in two parts, with the summation taken over $i$ from $1$ to $K$ in
the first part and over $i$ from $K+1$ to $M$ in the second, for
some fixed $K$, where $1 \le K \le M$. The idea behind this
splitting is to utilize the marginal distributions of the null
$p$-values from the first part through the inequality (\ref{equ5}), as we did
before, and the pairwise joint distributions of these $p$-values
from the second part through the following new inequality (to be
proved in the \hyperref[sec6]{Appendix}):
%
\begin{eqnarray}\label{equ20}
I (\widehat P_{(i)} \le t ) & \le& \frac{1}{i(i-1)} \sum
_{j=1}^{n_0} \sum_{j^{\prime}(\neq j)=1}^{n_0}I
\bigl(\max\{\widehat P_j, \widehat P_{j^{\prime}} \} \le t \bigr),
\end{eqnarray}
where $0< t<1$ is fixed, for all $i$ such that $2 \le i \le n_0$,
%
\begin{eqnarray}\label{equ21}
& & I \bigl( V > \max[\gamma R, k-1] \bigr)
\nonumber
\\
&&\qquad \le \sum_{i = 1 }^{K} \sum
_{j=1}^{n_0} \frac{1}{i \vee
k } I ( \widehat P_j \le \alpha_{{\bar m}(i)} ) I \bigl(\bigl\lfloor \gamma S/(1-
\gamma) \bigr\rfloor+1 = i \bigr)
\nonumber\\[-8pt]\\[-8pt]
&&\quad\qquad{} + \sum_{i = K +1 }^{M}
\sum_{j=1}^{n_0} \sum
_{l (\neq j)=1}^{n_0}
\frac{1}{(i \vee k)(i \vee k - 1)}I \bigl( \max (\widehat P_j, \widehat P_l ) \le \alpha_{\bar{m}(i)} \bigr)\nonumber
\\
&&\hspace*{122pt}{}\times   I \bigl(\bigl\lfloor\gamma
S/(1-\gamma) \bigr\rfloor+1 = i \bigr).\nonumber
\end{eqnarray}

Now, for each $j=1, \ldots, n_0$, the summation over $i$ in the
double-summation in~(\ref{equ21}) is equal to
%
\begin{eqnarray}\label{equ22}
& & \sum_{i = 1 }^{K} \biggl[
\frac{I  (
\widehat P_j \le\alpha_{{\bar m} (i)}  )}{i \vee k } - \frac{I
( \widehat P_j \le\alpha_{{\bar m}(i-1)}  )}{(i - 1) \vee k} \biggr] I \bigl(\bigl\lfloor\gamma S/(1-
\gamma) \bigr\rfloor+1 \ge i \bigr)
\nonumber
\\
&&\quad {} -  \frac{I
( \widehat P_j \le\alpha_{{\bar m}(K)}  )}{K \vee k} I \bigl(\bigl\lfloor\gamma S/(1-\gamma) \bigr
\rfloor+1 \ge K+1 \bigr)
\nonumber\\[-8pt]\\[-8pt]
&&\qquad \le \sum_{i = 1 }^{K}
\frac{I  ( \alpha_{{\bar m}(i-1)}
< \widehat P_j \le\alpha_{{\bar m}(i)}  )}{i \vee k}
\nonumber
\\
&&\quad\qquad{} -  \frac{I  ( \widehat P_j \le\alpha_{{\bar m}(K)}
)}{K \vee k} I \bigl(\bigl\lfloor\gamma S/(1-\gamma) \bigr
\rfloor+1 \ge K+1, M \ge K +1 \bigr)\nonumber
\end{eqnarray}
with $I(\widehat P_j \le\alpha_{{\bar m}(0)})/0 \vee k =0$, and
similarly for each $j \neq l$, the summation over $i$ in the
triple-summation in (\ref{equ21}) is less than or equal to
%
\begin{eqnarray}\label{equ23}
& & \sum_{i = K +2 }^{M} \biggl[
\frac{I  ( \max(\widehat P_j, \widehat P_l ) \le
\alpha_{\bar{m}(i)}  )}{(i \vee k)(i \vee k - 1)} - \frac{I  (
\max(\widehat P_j, \widehat P_l ) \le\alpha_{\bar{m}(i-1)}
)}{((i -1) \vee k)((i-1) \vee k- 1)} \biggr]
\nonumber
\\
&&\hspace*{26pt}{}\times I \bigl(\bigl\lfloor\gamma S/(1-\gamma) \bigr\rfloor+ 1 \ge i \bigr)
\nonumber
\\
&&\quad{} +  \frac{I  ( \max(\widehat P_j, \widehat P_l )
\le\alpha_{\bar{m}(K+1)}  )}{((K + 1) \vee k)[(K+1) \vee k -1]}\nonumber
\\
&&\hspace*{22pt}{}\times I \bigl(\bigl\lfloor\gamma S/(1-\gamma) \bigr
\rfloor+1 \ge K + 1, M \ge K +1 \bigr)
\\
&&\qquad \le \sum_{i = K +2 }^{M} \frac{I  ( \alpha_{\bar{m}(i-1)} < \max(\widehat P_j, \widehat P_l ) \le\alpha_{\bar{m}(i)}  ) }{(i \vee k)(i \vee k
- 1)}\nonumber
\\
&&\quad\qquad{} +  \frac{I  ( \max(\widehat P_j, \widehat P_l ) \le
\alpha_{\bar{m}(K+1)}  )}{((K + 1) \vee k)[(K+1) \vee k -1]}\nonumber
\\
&&\hspace*{44pt}{}\times  I \bigl(\bigl\lfloor\gamma S/(1-\gamma) \bigr
\rfloor+1 \ge K+1, M \ge K +1 \bigr).
\nonumber
\end{eqnarray}
In addition, by simple algebraic calculation, we have
%
\begin{eqnarray}\label{equ24}
& & \frac{I(\widehat{P}_j \vee\widehat{P}_l \le\alpha_{\bar{m}(K+1)}
)}{( (K+1) \vee k)( (K+1) \vee k-1)} - \frac{I(\widehat{P}_j \le
\alpha_{\bar{m}(K)} ) }{ (K \vee k)(n_0-1)}
\nonumber
\\
&&\qquad \le \frac{ (n_0-(K+1) \vee k  )I(\widehat{P}_j \le
\alpha_{\bar{m}(K)},\widehat{P}_l \le\alpha_{ \bar{m}(K+1)}
)}{((K+1) \vee k) ( (K+1) \vee k -1)(n_0-1)}
\\
&&\quad\qquad{} + \frac{I(\alpha_{\bar{m}(K)} < \widehat{P}_j \le\alpha_{\bar{m}(K+1)},
\widehat{P}_l \le\alpha_{\bar{m}(K+1)})}{( (K+1) \vee k)( (K+1) \vee k-1)}.\nonumber
\end{eqnarray}

Applying (\ref{equ22})--(\ref{equ24}) to (\ref{equ21}) and taking expectations of both
sides in (\ref{equ21}), we get
%
\begin{eqnarray}\label{equ25}
\gamma\mbox{-kFDP} & \le& \sum_{i = 1 }^{K}
\frac{ n_0(\alpha_{{\bar m}(i)} - \alpha_{{\bar m}(i-1)}  )}{i \vee k}
\nonumber
\\
&&{} + \sum_{i = K +2 }^{M} \frac{ n_0 (n_0 - 1) [ F
(\alpha_{\bar{m}(i)}, \alpha_{\bar{m}(i)}  )
- F  (\alpha_{\bar{m}(i-1)}, \alpha_{\bar{m}(i-1)}  )
]}{(i \vee k) (i \vee k - 1)}
\nonumber\\[-8pt]\\[-8pt]
& &{} + \frac{ n_0 (n_0 - 1) F  (\alpha_{\bar{m}(K+1)},
\alpha_{\bar{m}(K+1)}  ) }{((K + 1) \vee k) ((K+1) \vee k
-1)}I (M \ge K + 1 )
\nonumber
\\
& &{} - \frac{ n_0 F  (\alpha_{\bar{m}(K)},
\alpha_{\bar{m}(K+1)}  ) }{(K + 1) \vee k}I (M \ge K + 1 ).\nonumber
\end{eqnarray}

This\vspace*{1pt} inequality produces the next theorem, one of our main results in
this subsection, with $C_{n,\mathrm{SD}}^{(3)}(\beta)$ in that theorem being
defined as follows:
\begin{eqnarray*}
&& C_{k,n,\mathrm{SD}}^{(3)}(\beta)
\\
&&\qquad = \max_{ k \le n_0 \le n} \min
_{1 \le K
\le M} \Biggl\{ \sum_{i = 1 }^{K}
\frac{ n_0  [\alpha_{\bar
{m}(i)}^{\prime}(\beta) - \alpha_{\bar{m}(i-1)}^{\prime}(\beta)
]}{i \vee k}
\\
&&\hspace*{101pt}{} + \sum_{i = K +2 }^{M} \bigl( n_0(n_0-1)  \bigl[ F
\bigl(\alpha_{\bar{m}(i)}^{\prime}(\beta),
\alpha_{\bar{m}(i)}^{\prime}(\beta)  \bigr)
\\
&&\hspace*{196pt}{}- F  \bigl(\alpha_{\bar{m}(i-1)}^{\prime}(\beta), \alpha_{\bar
{m}(i-1)}^{\prime}(\beta)  \bigr)  \bigr]\bigr)
\\
&&\hspace*{143pt} /\bigl({(i \vee k) (i \vee k - 1)}\bigr)
\\
&&\hspace*{101pt}{} + \frac{ n_0(n_0 - 1) F
(\alpha_{\bar{m}(K+1)}^{\prime}(\beta),
\alpha_{\bar{m}(K+1)}^{\prime}(\beta)  ) }{((K + 1) \vee k)
((K+1) \vee k -1)}
\\
&&\hspace*{112pt}{}\times  I (M \ge K + 1 )
\\
&&\hspace*{131pt}{} -   \frac{ n_0 F
(\alpha_{\bar{m}(K)}^{\prime}(\beta),
\alpha_{\bar{m}(K+1)}^{\prime}(\beta)  ) }{(K + 1) \vee k}I (M \ge K + 1 ) \Biggr\},
\end{eqnarray*}
given a sequence of constants $0 = \alpha_{0}^{\prime}(\beta) \le\alpha
_1^{\prime}(\beta) \le\cdots\le\alpha_n^{\prime}(\beta)$ with a
fixed \mbox{$\beta\in(0,1)$}.
%
\begin{theorem}\label{th3.7}
Given any sequence of
critical constants $0 = \alpha_{0}^{\prime}(\beta) \le
\alpha_1^{\prime}(\beta) \le\cdots\le\alpha_n^{\prime}(\beta)$,
for a fixed $\beta\in(0,1)$, the stepdown procedure with the
critical values $\alpha_{i \vee k}, i=1, \ldots, n$, satisfying $\alpha
_{i \vee k} =
\alpha_{i \vee k}^{\prime}(\beta_{\mathrm{SD}}^*)$ and\break  $C_{k,n,\mathrm{SD}}^{(3)}(\beta
_{\mathrm{SD}}^*) =\alpha$, controls the $\gamma$-kFDP at $\alpha$
under Assumptions~\ref{ass1} and~\ref{ass3}.
\end{theorem}

We now derive a stepup analog of Theorem~\ref{th3.7} starting from the
following inequality, which is obtained from Lemma~\ref{le3.2} by splitting the
right-hand sum in the second inequality of that lemma into two parts,
as before, for a fixed $1 \le k \le K \le n_0$:
\begin{eqnarray*}
& & I\bigl(V > \max[\gamma R, k-1]\bigr)
\\
&&\qquad \le \sum_{j=1}^{n_0} \frac{I  (\widehat P_j
\le\alpha_{\tilde m(k-1)}  )}{k} +
\sum_{j=1}^{n_0} \sum
_{i=k}^{K} \frac{I  (\alpha_{\tilde m(i-1)}< \widehat P_j
\le\alpha_{\tilde m(i)}
)}{i}
\\
&&\quad\qquad{}+ \sum_{j=1}^{n_0} \sum
_{i=K + 1}^{n_0} \frac{I  ( \widehat R_2 \ge i ) I  (\alpha_{\tilde m(i-1)}< \widehat P_j
\le\alpha_{\tilde m(i)}  ) }{i}.
\end{eqnarray*}
Again, the idea behind this splitting is to capture the pairwise joint
distributions of the null $p$-values from the second part, and for
that, we use the following inequality, which can be seen to follow from
Lemma~\ref{le3.2} (see Remark~\ref{re3.2}):
%
\begin{eqnarray}\label{equ26}
& & I (\widehat R_{2} \ge r ) \le \sum_{l=1}^{n_0}
\Biggl(\frac{I  (\widehat P_l
\le\alpha_{\tilde m(r)}  )}{r} + \sum_{s=r+1}^{n_0}
\frac{I
(\alpha_{\tilde m(s-1)} < \widehat P_l \le\alpha_{\tilde m(s)}  )}{s} \Biggr).
\end{eqnarray}

Thus, we get
%
\begin{eqnarray}\label{equ27}
\hspace*{-5pt}&&  I \bigl( V > \max[\gamma R, k-1] \bigr)
\nonumber
\\
\hspace*{-5pt}&&\qquad  \le\sum_{j=1}^{n_0} \frac{I  (\widehat P_j \le\alpha_{\tilde m(k-1)}  )}{k}\nonumber
\\
\hspace*{-5pt}&&\!\! \quad\qquad{}  + \sum_{j=1}^{n_0} \sum
_{r=k}^{K} \frac{I
(\alpha_{\tilde m(r-1)}< \widehat P_j \le\alpha_{\tilde m(r)}  )}{r}
\nonumber\\[-8pt]\hspace*{-5pt} \\[-8pt]
\hspace*{-5pt}&&\!\!\quad\qquad{} + \sum_{j=1}^{n_0} \sum
_{r=K + 1}^{n_0} \frac{I(\alpha_{\tilde m(r-1)}< \widehat P_j \le\alpha_{\tilde m(r)})}{r^2}\nonumber
\\
\hspace*{-5pt}&&\!\!\quad\qquad{} + \sum_{j=1}^{n_0} \sum
_{l(\neq j)=1}^{n_0} \sum_{r= K +1}^{n_0}
\sum_{s=r+1}^{n_0} \frac{I  (\alpha_{\tilde m(r-1)}< \widehat P_j
\le\alpha_{\tilde m(r)}, \alpha_{\tilde m(s-1)}< \widehat P_l \le
\alpha_{\tilde m(s)}  )}{rs}\nonumber
\\
\hspace*{-5pt}&&\!\!\quad\qquad{} + \sum_{j=1}^{n_0} \sum
_{l(\neq j)=1}^{n_0} \sum_{r= K +1}^{n_0}
\!\frac{I  (\alpha_{\tilde m(r-1)}< \widehat P_j \le\alpha
_{\tilde m(r)}, \widehat P_l \le\alpha_{\tilde m(r)}  )}{r^2}.\nonumber
\end{eqnarray}
Taking expectations of both sides in (\ref{equ27}), we finally have
%
\begin{eqnarray}\label{equ28}
\gamma\mbox{-kFDP} &\le&\frac{n_0 \alpha_{\tilde m(k-1)}}{k}\nonumber
\\
&&{}  + \sum
_{r=k}^K \frac{n_0 (\alpha_{\tilde m(r)} - \alpha
_{\tilde m(r-1)} )} {r}
\nonumber
\\
&&{} + \sum_{r= K + 1}^{n_0} \frac{ n_0  (\alpha_{\tilde m
(r)} -\alpha_{\tilde m(r-1)}  )
}{r^2}
\\
&&{} +
\sum_{r= K + 1}^{n_0} \sum
_{s=r+1}^{n_0} \frac{ n_0  (
n_0 - 1  ) G  (\alpha_{\tilde m(r)}, \alpha_{\tilde m(s)}
) }{rs}
\nonumber
\\
&&{} + \sum_{r= K + 1}^{n_0} \frac{ n_0  ( n_0 - 1  )
( F(\alpha_{\tilde m (r)}, \alpha_{\tilde m (r)} ) -
F(\alpha_{\tilde m (r)}, \alpha_{\tilde m (r-1)}) ) }{r^2},\nonumber
\end{eqnarray}
where
\begin{eqnarray}
G(\alpha_{r}, \alpha_{s}) & = & F (\alpha_r,
\alpha_s ) - F (\alpha_{r-1}, \alpha_s ) - F
(\alpha_r, \alpha _{s-1} ) + F (\alpha_{r-1},
\alpha_{s-1} ).
\nonumber
\end{eqnarray}

Our next main result of this subsection follows from the inequality
(\ref{equ28}), with $C_{n,\mathrm{SU}}^{(3)}(\beta)$ in that result being defined as follows:
\begin{eqnarray*}
\hspace*{-4.5pt}&& C_{k,n,\mathrm{SU}}^{(3)}(\beta)
\\[-4pt]
\hspace*{-4.5pt}&&\qquad = \max_{ k \le n_0 \le n} \min
_{k \le K \le n_0 } \Biggl\{\frac{n_0 \alpha_{\tilde m(k-1)}^{\prime
}(\beta)}{k} + \sum
_{r=k}^{K} \frac{n_0  [\alpha_{\tilde m
(r)}^{\prime}(\beta) - \alpha_{\tilde m(r-1)}^{\prime}(\beta)
] }{r}
\\[-2pt]
\hspace*{-4.5pt}&&\hspace*{102pt}{}   + \sum_{r= K + 1}^{n_0} \Biggl(
\frac{ n_0[\alpha_{\tilde m (r)}^{\prime}(\beta) - \alpha_{\tilde m(r-1)}^{\prime}(\beta)  ]}{r^2} \nonumber
\\[-2pt]
\hspace*{-4.5pt}&&\hspace*{149pt} {}+ \sum_{s=r+1}^{n_0}
\frac{n_0(n_0-1) G  (\alpha_{\tilde m (r)}^{\prime}(\beta),
\alpha_{\tilde m (s)}^{\prime}(\beta)
) }{rs}
\\[-2pt]
\hspace*{-4.5pt}&&\hspace*{149pt}{} +  \bigl( n_0  ( n_0 - 1  )
\\[-2pt]
&&\hspace*{162pt}{}\times  \bigl[
F\bigl(\alpha_{\tilde m (r)}^{\prime}(\beta), \alpha_{\tilde m
(r)}^{\prime}(\beta) \bigr)
\\[-4pt]
\hspace*{-4.5pt}&&\hspace*{192pt}{} - F\bigl(\alpha_{\tilde m (r)}^{\prime}(\beta), \alpha_{\tilde m (r-1)}^{\prime}(\beta
)\bigr) \bigr] \bigr)
/ {r^2} \Biggr) \Biggr\},
\end{eqnarray*}
for any given sequence of constants $0 = \alpha_{0}^{\prime}(\beta) \le
\alpha_1^{\prime}(\beta) \le\cdots\le\alpha_n^{\prime}(\beta)$, for
a fixed $\beta\in(0,1)$.\vspace*{-1pt}

\begin{theorem}\label{th3.8}
Given any sequence of critical
constants $0 = \alpha_{0}^{\prime}(\beta) \le
\alpha_1^{\prime}(\beta) \le\cdots\le\alpha_n^{\prime}(\beta)$,
for a fixed $\beta\in(0,1)$, the stepup procedure with the
critical values $\alpha_{i\vee k}, i=1, \ldots, n$, satisfying $\alpha
_i =
\alpha_i^{\prime}(\beta_{\mathrm{SU}}^*)$ and $C_{k,n,\mathrm{SU}}^{(3)}(\beta_{\mathrm{SU}}^*) =
\alpha$, controls the $\gamma$-kFDP at $\alpha$ under Assumptions~\ref{ass1} and~\ref{ass3}.
\end{theorem}
%

%
\begin{remark}\label{re3.5}
\citeauthor{r16} proved the following two results in (\citeyear{r16}) and
(\citeyear{r17}), respectively, based on marginal $p$-values under arbitrary
dependence: the $\gamma$-FDP of the stepdown procedure
with critical values $\alpha_i, i=1, \ldots, n$, satisfies\vspace*{-3pt}
%
\begin{eqnarray}\label{equ29}
\gamma\mbox{-FDP} & \le& \max_{1 \le n_0 \le n } \Biggl\{
n_0 \sum_{i =1}^{M}
\frac{\alpha_{\bar{m}(i)} -
\alpha_{\bar{m}(i - 1)}}{i} \Biggr\};
\end{eqnarray}
whereas
the $\gamma$-FDP of the stepup procedure with the same set of
critical values satisfies\vspace*{-3pt}
%
\begin{eqnarray}\label{equ30}
\gamma\mbox{-FDP} & \le& \max_{1 \le n_0 \le n } \Biggl\{
n_0 \sum_{i=1}^{n_0}
\frac{\alpha_{\tilde m(i)} - \alpha
_{\tilde m(i-1)}}{i} \Biggr\}.
\end{eqnarray}
These upper bounds are always larger than the
corresponding upper bounds of the $\gamma$-FDP we derive here, as
seen by letting $k=1, K=M$ in (\ref{equ25}) and $k=1, K=n_0$ in (\ref{equ28}),
respectively. Thus, theoretically, the stepdown and stepup
$\gamma$-FDP controlling procedures introduced in Theorems~\ref{th3.7} \mbox{(with
$k=1$)}\break  and~\ref{th3.8} (with $k=1$), respectively, are always more powerful
than the corresponding ones given by \citeauthor{r16} in (\citeyear{r16})
and (\citeyear{r17}), respectively.\vadjust{\goodbreak}
\end{remark}

\section{Simulation studies}\label{sec4}
We ran extensive simulations numerically examining the performances
of different procedures proposed in the above section in comparison
with their relevant competitors under different settings for
the\vadjust{\goodbreak}
parameters, $\pi_0$, $\gamma$, $k$ and the strength of positive
dependence, and having considered all or just one of three special types
of positive dependence structure---uniform pairwise dependence,
clumpy dependence and autoregressive of order one [AR(1)]
dependence. The results were graphically summarized in twelve
figures, and the main findings in those graphs are described in this
section. However, we present here the figures that pertain to the
uniform pairwise dependence, while the rest are presented, for lack
of space here, in the supplementary material [\citet{r71}].

Note that, except in the procedures in Theorems~\ref{th3.1} and~\ref{th3.4}, which
have been developed directly from the Lehmann--Romano (LR) critical
values $\alpha'_i(\beta) = \frac{(\lfloor\gamma i \rfloor+ 1)\beta} {
n + \lfloor\gamma i \rfloor+ 1 - i}, i = 1, \ldots, n$, the critical\vspace*{1pt}
values in all other procedures can be chosen
arbitrarily before being rescaled appropriately to ensure a control
over the $\gamma$-FDP or $\gamma$-kFDP. In many of our simulations,
we had chosen the same LR critical values $\alpha'_i(\beta)$ in these
other procedures
with $\beta$ being rescaled according to the formulas in the
corresponding theorems. We will refer to a procedure, except the
stepwise one in Theorem~\ref{th3.1}, as simply LR-type procedure whenever it
is directly or indirectly based on the LR critical values.
Similarly, by BH- and GBS-type stepwise procedures that we will use
in some simulations, we mean that the critical values of the
procedure in that procedure are obtained by rescaling the original
BH or GBS critical values according to the formula given in the
corresponding theorem.

A part of our simulation study was geared toward answering the following
two questions:
\begin{longlist}[(Q1)]
\item[(Q1)] When controlling the $\gamma$-FDP assuming positive dependence,
how good is the improvement supposedly offered by the newly proposed
LR stepup procedure in Theorem~\ref{th3.1} over the original LR stepdown
procedure?
\item[(Q2)] When controlling the $\gamma$-FDP assuming arbitrary dependence,
how do the newly suggested LR-type stepdown and stepup procedures in
Theorems~\ref{th3.7} and~\ref{th3.8}, respectively, with $k=1$, incorporating
pairwise correlation information perform compared to the
corresponding existing LR-type stepdown and stepup procedures in
\citeauthor{r16} (\citeyear{r16,r17}) that do not incorporate such pairwise
correlation information?
\end{longlist}

The performance of each procedure is judged, while answering (Q1) and
(Q2), in terms of how well the $\gamma$-FDP is controlled at the
desired level and also the average power, which is the expected
proportion of false nulls that are rejected, under varying $\pi_0$,
$\gamma$, and the strength of positive dependence.

To simulate the values of $\gamma$-FDP and average power for each of
the methods referred to in (Q1) and (Q2), we first generated $n$
dependent normal random variables $N(\mu_i, 1), i = 1, \ldots, n$,
with $\pi_0 n$ of the $\mu_i$'s being equal to $0$ and the rest
being equal to $d =\sqrt{10}$, and a correlation matrix $\Gamma$.
The following three different types of $\Gamma$ were considered for
(Q1): (i) $\Gamma= (1-\rho) I_n+ \rho1_n 1_n^{\prime}$, in case of
uniform\vadjust{\goodbreak} pairwise dependence, (ii) $\gamma= I_{{n}/{s}} \otimes
[ (1-\rho)I_s+ \rho1_s 1_s^{\prime}  ]$, in case of
block dependence with the block size $s$, and (iii) $\Gamma=
((\rho^{|i-j|}))$, in case of AR(1) dependence, where $1_n=(1,
\ldots, 1)^{\prime}$; whereas, for (Q2), the $\Gamma$ of the type (i)
was considered. In each case, $\rho$ was nonnegative. We then
applied each method to the generated data to test $H_i\dvtx  \mu_i = 0$
against $K_i\dvtx  \mu_i \neq0$ simultaneously for $i =1, \ldots, n$, at
level $\alpha= 0.05$. We repeated the above two steps 2000
times.

%
\begin{figure}

\includegraphics{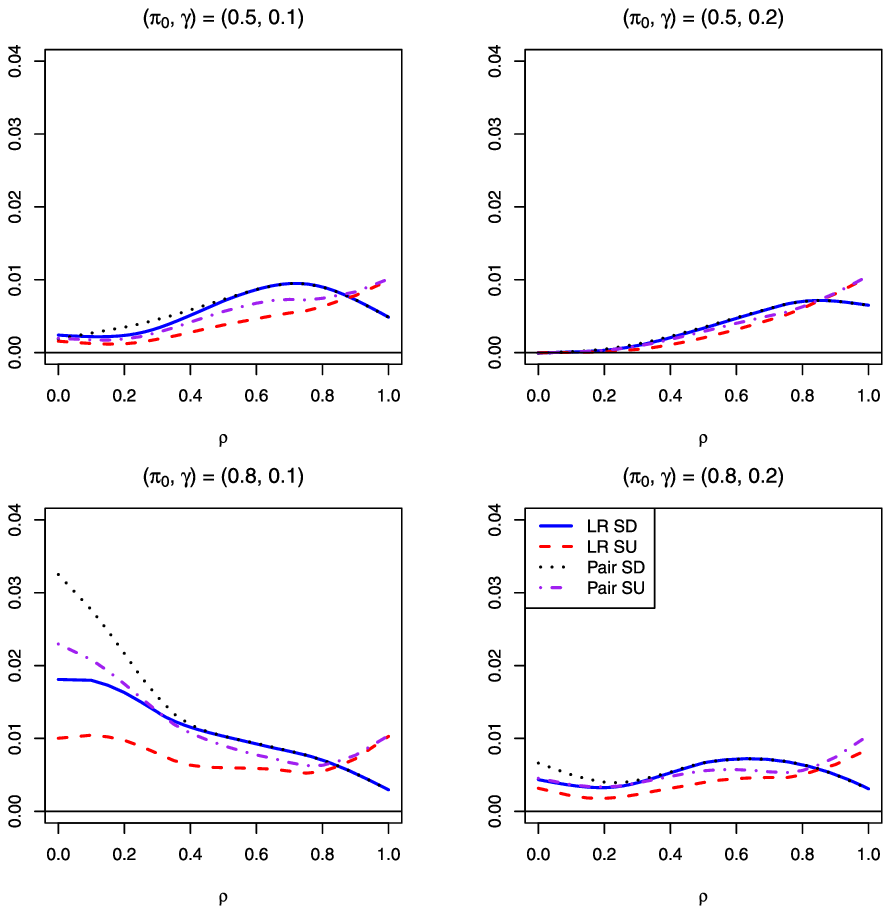}

\scriptsize{(a)~Simulated $\gamma$-FDP}

\caption{Simulated values of $\gamma$-FDP and average power of the existing
LR-type stepdown (LR SD) and stepup (LR SU) $\gamma$-FDP procedures
in Theorems~\protect\ref{th3.5} and~\protect\ref{th3.6} (with $k=1$) and the newly suggested
LR-type stepdown (Pair SD) and stepup (Pair SU) $\gamma$-FDP
procedures in Theorems \protect\ref{th3.7} and \protect\ref{th3.8} (with $k=1$), all developed assuming
arbitrary dependence, for $n=50$ and $\alpha= 0.05$.}\label{fig2}
\end{figure}

\setcounter{figure}{1}
\begin{figure}

\includegraphics{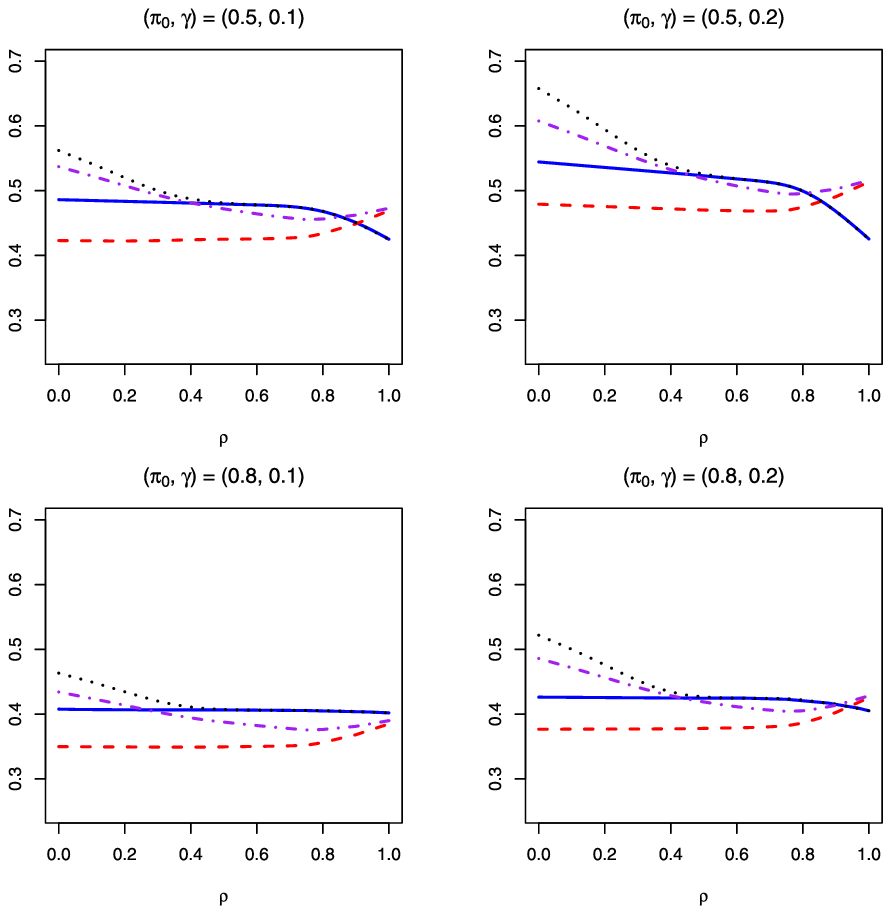}
\vspace*{-1pt}

\scriptsize{(b) Simulated average power}
\caption{(Continued).}\vspace*{-4pt}
\end{figure}

Figure~\ref{fig1} and Figures S.1--S.3 (in the supplementary material [\citet{r71}]) provide
an answer to (Q1) and Figure~\ref{fig2} answers (Q2). As seen from Figure~\ref{fig1}, when
the underlying
test statistics have a common positive correlation, the
newly introduced stepup $\gamma$-FDP procedure in Theorem~\ref{th3.1} improves
the power
of the original Lehmann--Romano stepdown
procedure. This improvement is quite noticeable when the correlation is
high. When the underlying
test statistics are block or AR(1) dependent, the stepup procedure,
as expected, does still have better power, as seen from Figures
S.1--S.3. However, in case of block dependence, the larger the block
size, the more significant seems to be the power improvement;
whereas, in case of AR(1) dependence, the power improvement seems to
be only significant when the dependence is high and the proportion
of true nulls is not large. In addition, as seen from Figure~\ref{fig1}, for
the original LR stepdown procedure and its stepup analogue, they
behave very differently when correlation $\rho$ is close to one,
which corroborates the observation of \citet{r19},
and their apparent worst performances in terms of the $\gamma$-FDP
control seem to depend on the values of $\pi_0$ and $\gamma$.

From Figure~\ref{fig2}, we see that when controlling the $\gamma$-FDP
assuming arbitrary dependence, the performances of the existing
LR-type Romano--Shaikh stepdown and stepup procedures can be
significantly improved by utilizing the pairwise correlation
information via the use of the newly suggested LR-type stepdown and
stepup $\gamma$-FDP procedures in Theorems~\ref{th3.7} and~\ref{th3.8},
respectively, with $k=1$, when the underlying test statistics are
slightly or moderately correlated with a common correlation.

Our next set of simulations was run with a view to investigating the
performances
of the proposed stepwise $\gamma$-kFDP controlling procedures in the
setting of a common pairwise positive dependence. Specifically, we
investigated the following two questions:
\begin{longlist}[(Q3)]
\item[(Q3)] When controlling the $\gamma$-kFDP assuming positive
dependence, how well the LR-type stepdown and stepup procedures in
Theorem~\ref{th3.4} incorporating\vadjust{\goodbreak} pairwise correlation information perform
compared to the LR-type stepdown and stepup procedures in Theorems~\ref{th3.2}
and~\ref{th3.3}, respectively, that do not incorporate such pairwise
correlation information?
\item[(Q4)] When controlling the $\gamma$-kFDP assuming arbitrary
dependence, how well do the LR-type stepdown and stepup procedures in
Theorems~\ref{th3.7} and~\ref{th3.8}, respectively, incorporating pairwise correlation
information perform compared to the LR-type stepdown and stepup
procedures in Theorems~\ref{th3.5} and~\ref{th3.6}, respectively, that do not
incorporate such pairwise correlation information?
\end{longlist}

The performance of each procedure is judged, while answering (Q3) and (Q4),
in terms of how well the $\gamma$-kFDP is controlled at the desired
level and also the average power under varying $\pi_0$, $k$ and the
strength of positive dependence. We used the simulation settings for (Q3)
and (Q4) that are same as in answering (Q1) and (Q2), respectively, but
considering only the equi-correlated normal case.

%
\begin{figure}

\includegraphics{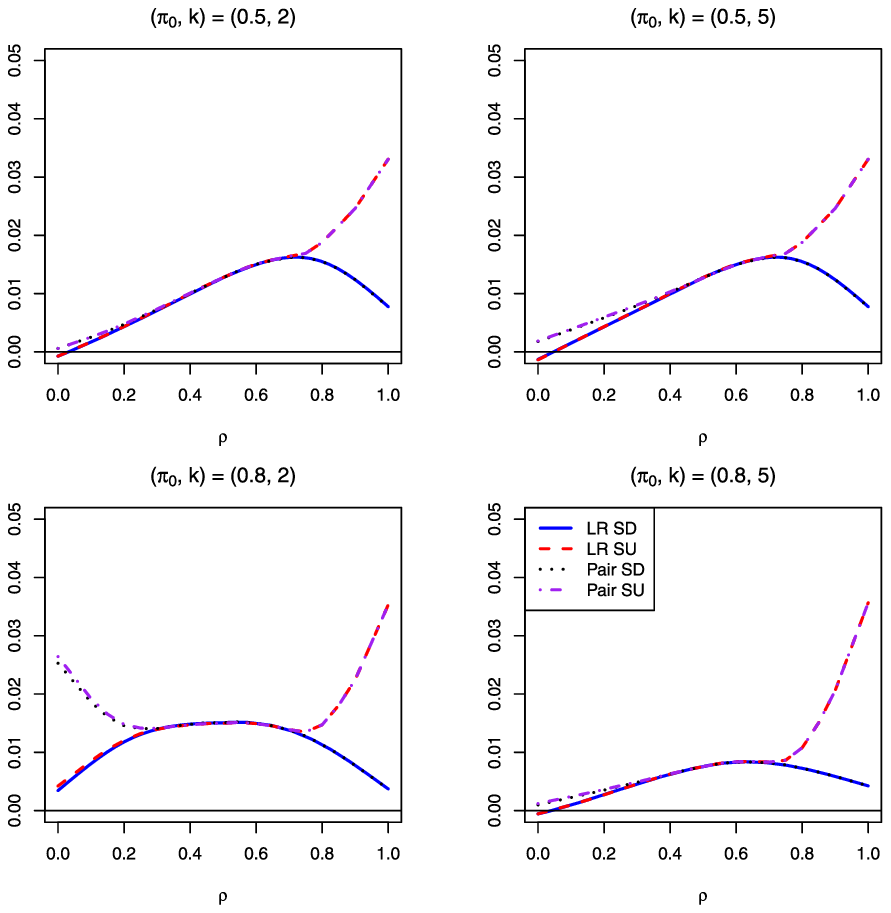}

\scriptsize{(a)~Simulated $\gamma$-kFDP}\vspace*{-3pt}
\caption{Simulated values of $\gamma$-kFDP and average power of the
LR stepdown (LR SD) and stepup (LR SU)
$\gamma$-kFDP procedures in Theorems~\protect\ref{th3.2}
and~\protect\ref{th3.3} and the LR-type stepdown
(Pair SD) and stepup (Pair SU) $\gamma$-kFDP procedures in Theorem \protect\ref{th3.4}, all
developed assuming positive dependence, for $n=100, \gamma= 0.1$ and
$\alpha= 0.05$.}\label{fig3}\vspace*{-3pt}
\end{figure}

\setcounter{figure}{2}
\begin{figure}

\includegraphics{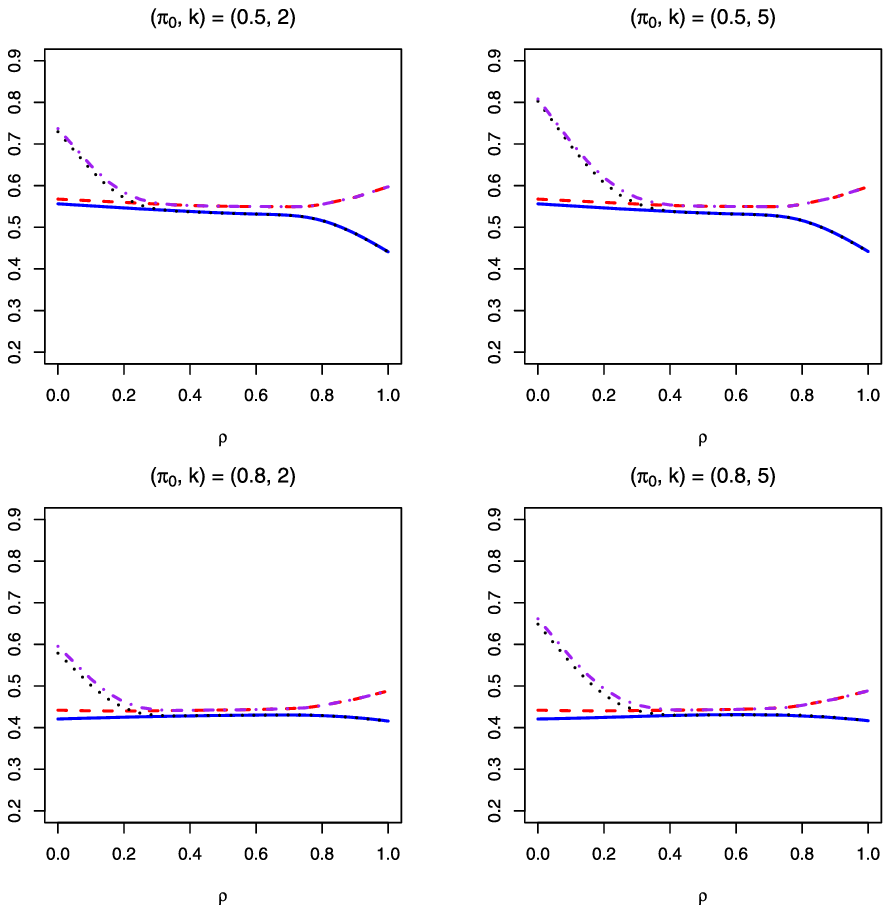}

\scriptsize{(b) Simulated average power}
\caption{(Continued).}
\end{figure}

%
\begin{figure}

\includegraphics{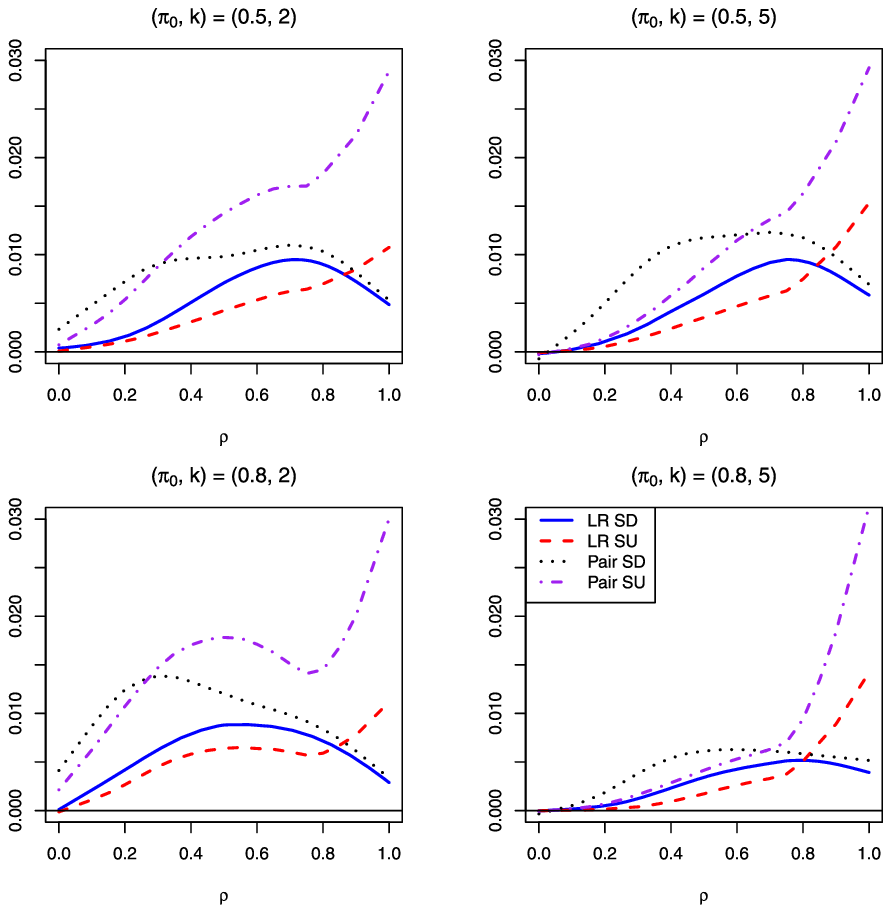}

\scriptsize{(a)~Simulated $\gamma$-kFDP}
\caption{Simulated values of $\gamma$-kFDP and average power of the LR-type
stepdown (LR SD) and stepup (LR SU)
$\gamma$-kFDP procedures in Theorems~\protect\ref{th3.5} and~\protect\ref{th3.6}
and the LR-type stepdown (Pair SD)
and stepup (Pair SU) $\gamma$-kFDP procedures in Theorems \protect\ref{th3.7} and \protect\ref{th3.8}, all
developed assuming arbitrary dependence, for $n=50, \gamma= 0.1$ and
$\alpha= 0.05$.}\label{fig4}
\end{figure}

\setcounter{figure}{3}
%
\begin{figure}

\includegraphics{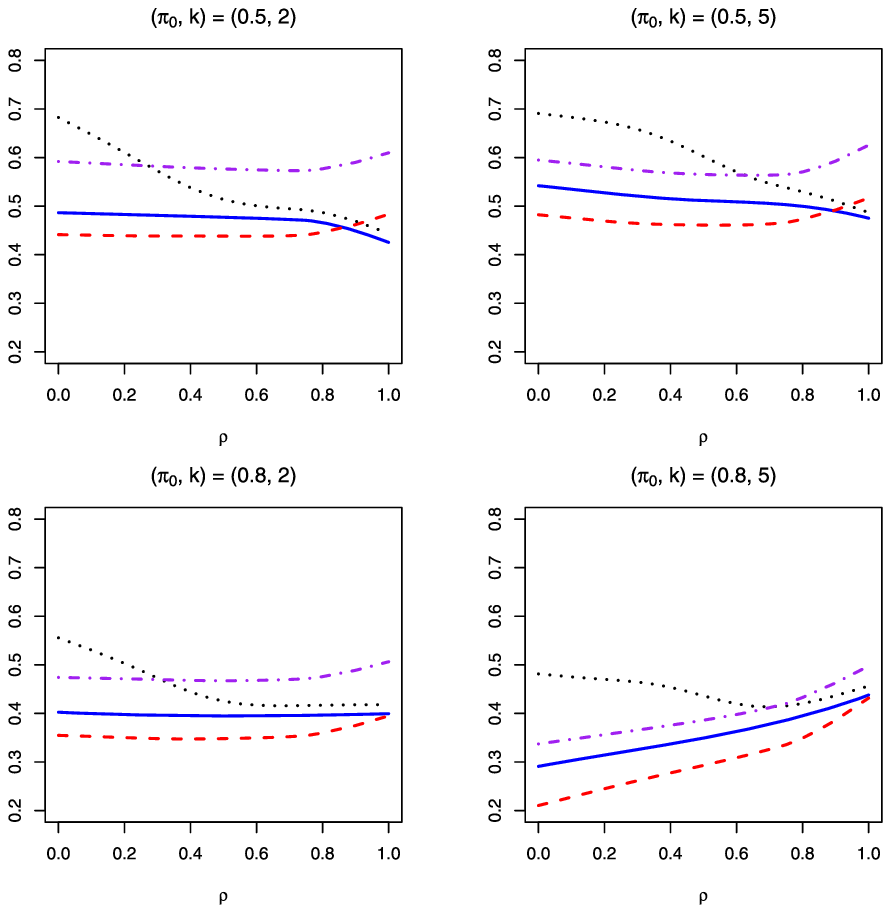}

\scriptsize{(b) Simulated average power}
\caption{(Continued).}
\end{figure}

Figures~\ref{fig3} and~\ref{fig4} provide answers to (Q3) and (Q4), respectively. From Figure~\ref{fig3},
we see that when controlling the $\gamma$-kFDP assuming positive
dependence, the stepwise $\gamma$-kFDP procedure in Theorems~\ref{th3.2} or~\ref{th3.3},
which is based only on the marginal $p$-values, seem to perform well,
but it can be significantly improved by utilizing the pairwise
correlation information via the use of the corresponding stepwise
procedure in Theorem~\ref{th3.4} when the underlying test statistics are weakly
correlated. However, when the test statistics are
strongly correlated, this stepwise procedure incorporating such
pairwise correlations has almost the same power performance as the
corresponding stepwise procedure based only on the marginal \mbox{$p$-}values.
Of course, such phenomenon has been noted before in the context of
other generalized error rates [\citet{r34}]. Figure~\ref{fig4},
however, reveals an interesting picture. It seems to say that when
controlling the $\gamma$-kFDP assuming arbitrary dependence, the
LR-type stepwise procedure in Theorems~\ref{th3.5} or~\ref{th3.6} based only on the
marginal $p$-values can be made consistently more powerful by utilizing
the pairwise correlation information through the use of the
corresponding LR-type stepwise procedure in Theorems~\ref{th3.7} or~\ref{th3.8}, with
the power gaps still being quite significant even when the test
statistics are highly correlated.

Looking at all the these seven figures, it becomes clear that given a
choice of~$\gamma$, the performance of an LR-type stepwise procedure,
particulary in terms of controlling the $\gamma$-FDP or $\gamma$-kFDP,
is affected not only by dependence but also by~$\pi_0$.

We also did some simulations to examine the following question:
\begin{longlist}[(Q5)]
\item[(Q5)] How do the newly suggested BH- and GBS-type $\gamma$-FDP
stepup procedures assuming positive dependence in Theorem~\ref{th3.3}\vadjust{\goodbreak} with $k
= 1$ perform compared to the corresponding BH- and GBS-type $\gamma$-FDP
stepdown procedures obtained from Theorem~\ref{th3.2}?
\end{longlist}

We used the same simulation settings involving three different types
of positive dependence structure as in answering (Q1). From Figures
S.4--S.7 (in the supplementary material [\citet{r71}]) that answers (Q5), we see that
the BH- or GBS-type stepup and stepdown procedures have the similar
behaviors as the LR-type procedures. Generally, when the underlying
test statistics are highly correlated, the power improvements of the
stepup procedures over the corresponding stepdown procedures are
always quite significant. For other cases, the power improvement
depends on the dependence structure of the test statistics. In
addition, an interesting observation is that the BH-type stepwise
procedures are always more powerful than the corresponding GBS-type
procedures.

\setcounter{figure}{4}
%
\begin{figure}

\includegraphics{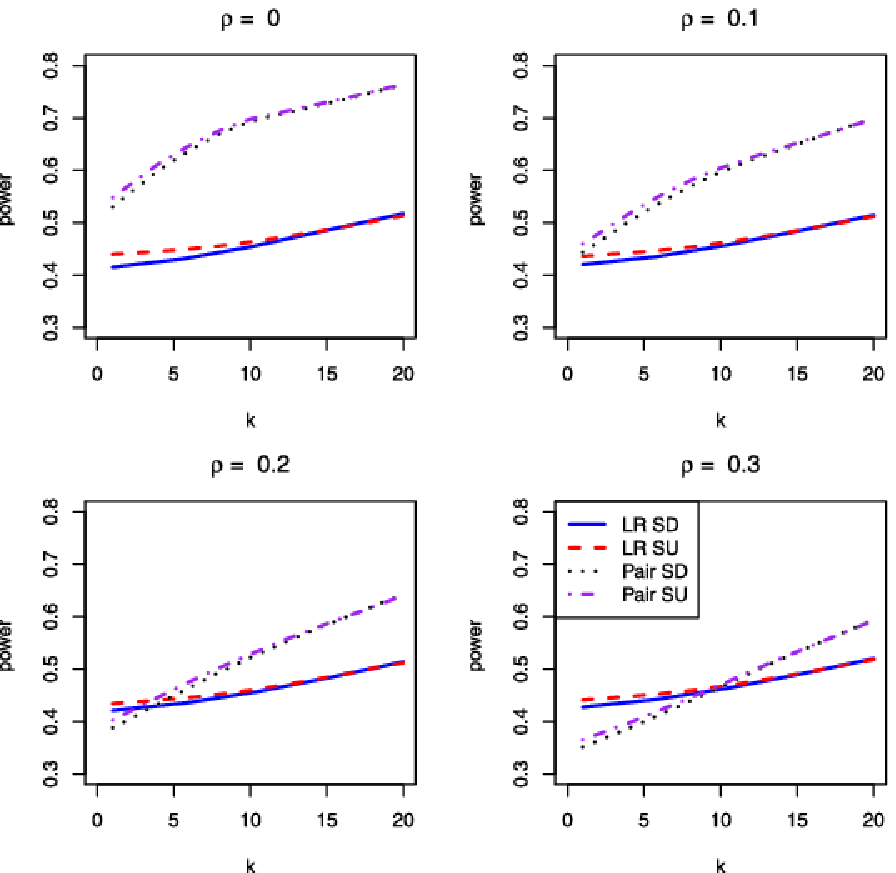}

\caption{Simulated average power of the LR stepdown (LR SD) and stepup
(LR SU)
$\gamma$-kFDP procedures in Theorems \protect\ref{th3.2} and \protect\ref{th3.3} and the LR-type
stepdown (Pair SD) and stepup (Pair SU) $\gamma$-kFDP procedures in
Theorem \protect\ref{th3.4} with respect to different values of $k$, all developed
assuming positive dependence, for $n = 100, \pi_0 = 0.8$, $\gamma
= 0.1$ and $\alpha= 0.05$.}\label{fig5}
\end{figure}

Our last set of simulations was carried out to investigate the following:
\begin{longlist}[(Q6)]
\item[(Q6)] As a $\gamma$-kFDP procedure under positive dependence, how does
the \mbox{LR-}type stepwise procedure in Theorem~\ref{th3.4}
incorporating\vadjust{\goodbreak}
pairwise correlation perform in terms of power with increasing $k$ and
strength of dependence,
compared to the corresponding LR-type stepwise
procedure in Theorems~\ref{th3.2} or~\ref{th3.3} that do not incorporate such pairwise
correlation
information?
\end{longlist}

We used the same simulation setting as in answering (Q3). From Figure~\ref{fig5}
that answers this question, we see
that the power of each of these LR-type stepwise \mbox{$\gamma$-}kFDP
procedures increases with $k$, as expected. The power gap between
the stepwise $\gamma$-kFDP procedure in Theorem~\ref{th3.4} and the
corresponding stepwise \mbox{$\gamma$-}kFDP procedure in Theorems~\ref{th3.2} or
\ref{th3.3} gets wider with increasing $k$. The~stepwise procedures in
Theorem~\ref{th3.4} are more powerful than the corresponding stepwise
procedures in Theorems~\ref{th3.2} and~\ref{th3.3}, irrespective of choice of $k$ if
the \mbox{underlying} test
statistics are weakly correlated and with properly chosen $k$ if
these statistics are moderately correlated.

\section{Concluding remarks}\label{sec5} The paper is motivated by the need to advance
the theory of FDP control which is still underdeveloped despite
being well accepted by the multiple testing research community. Our
focus has been two-fold: (i) enlarging the class of procedures
controlling the $\gamma$-FDP, the existing notion of FDP control,
and (ii) generalizing this notion to one that is often more
appropriate and powerful---with improving some of the currently
available results under certain dependence situations being the
overreaching goal. We have given a large class of procedures
controlling the $\gamma$-FDP and its generalization under different
dependence assumptions, and numerical evidences showing superior
performances of the proposed procedures compared to those they
intend to improve under some dependence cases, although these
proposed procedures themselves, like their competitors, are still quite
conservative.

There is scope of doing further research in the context of what we
discuss in this paper. We have defined the $\gamma$-kFDP, for the
first time in this paper, with the idea of introducing a more
powerful notion of error rate than the $\gamma$-FDP under
dependence. We have proposed several procedures controlling the
$\gamma$-kFDP and given numerical evidence of their power
superiority over the corresponding $\gamma$-FDP controlling
procedures for some specific values of $k$ and under certain
dependence situations. Although a deeper understanding of $\gamma$-kFDP
under dependence, particularly, how the choice of $k$
depends on correlations, would require studying distributional
properties of FDP or kFDP under dependence, an
area still less developed, we have provided some insight into it
through additional simulations whose findings are reported in the
supplementary material [\citet{r71}]. In particular, it has been noted that the
difference between controlling $\gamma$-kFDP and $\gamma$-FDP and the
stipulated power gain in using a $\gamma$-kFDP procedure over the
corresponding $\gamma$-FDP procedure may not be realized until $k/n$
reaches a certain critical point. Once this point is reached, the power
gain can be expected to steadily increase with $k/n$. Some idea about
the choice of $k$ relative to $n$ under different types and varying
strengths of dependence has also been provided.


\begin{appendix}\label{sec6}
\section*{Appendix}
\begin{pf*}{Proof of Lemma~\ref{le3.1}} First, note that
%
\begin{eqnarray}\label{equ31}
&& I \bigl( V > \max[\gamma R,k-1] \bigr)
\nonumber
\\
&&\qquad =  I \bigl(V > \max\bigl[ \gamma(V +S),k-1\bigr] \bigr)
\nonumber\\[-8pt]\\[-8pt]
&&\qquad =  I \bigl( V \ge\max\bigl\{\bigl\lfloor\gamma S/(1-\gamma) \bigr\rfloor +1,
k \bigr\} \bigr)
\nonumber
\\
&&\qquad =  \sum_{i = 1 }^{M} I \bigl( V \ge i \vee
k, \bigl\lfloor\gamma S/(1-\gamma) \bigr\rfloor+1 = i \bigr).\nonumber
\end{eqnarray}
Also, for a stepdown procedure with the critical
constants $\alpha_i$'s, we have
%
\begin{eqnarray}\label{equ32}
& & I \bigl( V \ge i \vee k, \bigl\lfloor\gamma S/(1-\gamma) \bigr\rfloor+1 = i
\bigr)\nonumber
\\
&&\qquad  =  I \bigl( R \ge i \vee k +S, \bigl\lfloor\gamma S/(1-\gamma) \bigr\rfloor+1
= i \bigr)\nonumber
\\
&&\qquad =  I \bigl( P_{(1)} \le\alpha_1, \ldots,
P_{(i \vee k +S)} \le\alpha_{i \vee k +S}, \bigl\lfloor\gamma S/(1-\gamma)
\bigr\rfloor+1 = i \bigr)
\nonumber\\[-8pt]\\[-8pt]
&&\qquad \le I \bigl( \widehat P_{(1)} \le\alpha_{1+S}, \ldots, \widehat P_{(i \vee k)} \le\alpha_{i \vee k +S}, \bigl\lfloor \gamma S/(1-\gamma)
\bigr\rfloor+1 = i \bigr)
\nonumber
\\
&&\qquad \le I \bigl( \widehat P_{(1)} \le\alpha_{1+m(i)}, \ldots, \widehat P_{(i \vee k)} \le \alpha_{i \vee k +m(i)}, \bigl\lfloor\gamma S/(1-\gamma) \bigr
\rfloor+1 = i \bigr)
\nonumber
\\
&&\qquad \le I \bigl( \widehat P_{(i \vee k )} \le\alpha_{i
\vee k +m(i)}, \bigl\lfloor
\gamma S/(1-\gamma) \bigr\rfloor+1 = i \bigr).\nonumber
\end{eqnarray}
Combining (\ref{equ31}) and (\ref{equ32}), we get the lemma.
\end{pf*}

\begin{pf*}{Proof of Lemma~\ref{le3.2}}
Since $ V \ge R-n_1$, we have
\begin{eqnarray*}
&& I \bigl(V > \max[\gamma R, k-1] \bigr)
\\
&&\qquad  = I \bigl(V \ge\max\bigl\{ \lfloor
\gamma R \rfloor+1, k \bigr\}, V \ge R-n_1 \bigr)
\\
&&\qquad =  I \Biggl( \bigcup_{j = k}^{n_0} \bigl\{
\widehat P_{(j)} \le \alpha_R, V =j,\lfloor \gamma R \rfloor+1
\le j, R \le j +n_1 \bigr\} \Biggr)
\\
&&\qquad \le I \Biggl( \bigcup_{j = k}^{n_0} \bigl\{
\widehat P_{(j)} \le \alpha_R, R \le\tilde m(j) \bigr\}
\Biggr)
\\
&&\qquad \le I \Biggl( \bigcup_{j = k}^{n_0} \{
\widehat P_{(j)} \le\alpha_{\tilde m(j)} \} \Biggr)
\\
&&\qquad =  \sum_{i = 1}^{n_0} \sum
_{j =
k}^{n_0} \frac{ I  ( \widehat P_i \le\alpha_{\tilde m(j)}, \widehat R_2 = j  ) }{
j}.
\end{eqnarray*}
This is the first inequality. The second inequality can be proved as follows:
\begin{eqnarray}
&& \sum_{i = 1}^{n_0} \sum
_{j = k}^{n_0} \frac{ I
( \widehat P_i \le\alpha_{\tilde m(j)}, \widehat R_2 = j  ) }{
j}
\nonumber
\\
&&\qquad = \sum_{i = 1}^{n_0} \sum
_{j = k}^{n_0} \frac{ I
( \widehat P_i \le\alpha_{\tilde m(j)}, \widehat R_2 \ge j  ) }{ j } - \sum
_{i = 1}^{n_0} \sum_{j =
k}^{n_0}
\frac{ I  ( \widehat P_i \le
\alpha_{\tilde m(j)}, \widehat R_2 \ge j+1  ) }{ j }
\nonumber
\\
&&\qquad \le  \sum_{i = 1}^{n_0} \frac{ I  ( \widehat P_i \le\alpha
_{\tilde m(k)}, \widehat R_2 \ge k  ) }{ k } +
\sum_{i = 1}^{n_0} \sum
_{j = k+1}^{n_0} \frac{ I  ( \alpha_{\tilde m(j-1)} < \widehat P_i \le\alpha_{\tilde m(j)}, \widehat R_2 \ge j  ) }{ j }.
\nonumber
\end{eqnarray}
Thus, the lemma is proved.
\end{pf*}

%
\begin{proposition}\label{pr6.1}
Let $M$ and $m(i)$, for $i=1, \ldots, M$, be defined as in Lemma~\ref{le3.1}
and $\tilde{m}(i)$ for $i=1, \ldots, n_0$ be defined as in Lemma~\ref{le3.2}. Then, for given set of critical constants,
\[
\alpha'_i = \frac{(\lfloor\gamma i \rfloor+ 1) \alpha}{ n +
\lfloor\gamma i \rfloor+ 1 - i},\qquad i=1, \ldots,
n,
\]
we\vspace*{2pt} have $C_{k, n,\mathrm{SD}}^{(1)} = C_{k, n, \mathrm{SU}}^{(1)} = \alpha$ when $k=1$,
where $C_{k, n,\mathrm{SD}}^{(1)}$ and $C_{k, n,\mathrm{SU}}^{(1)}$ are, respectively,
defined as in Theorems~\ref{th3.2} and~\ref{th3.3}.
\end{proposition}

\begin{pf}
We first prove that $C_{k, n,\mathrm{SD}}^{(1)} = \alpha$ when $k=1$.
From the definition of $m(i)$, we have
\[
i -1 \le\frac{\gamma m(i)}{1 - \gamma} < i.
\]
Thus,
\[
i-1 \le i - (1 - \gamma) \le\gamma\bigl(i + m(i)\bigr) < i.
\]
Hence,
%
\begin{equation}\label{equ33}
\bigl\lfloor\gamma\bigl(i + m(i)\bigr) \bigr\rfloor+ 1 = i.
\end{equation}
Based on (\ref{equ33}), we have
%
\begin{equation}\label{equ34}
\qquad \frac{n_0 \alpha'_{i + m(i)}}{i} = \frac{n_0 (\lfloor\gamma(i + m(i)) \rfloor+ 1)\alpha}{i(n + \lfloor
\gamma(i + m(i)) \rfloor+ 1
- i - m(i))} = \frac{n_0 \alpha}{n - m(i)} \le\alpha.
\end{equation}
Here, the inequality follows from the facts that $m(i) \le n_1$ and
$n_0 + n_1 = n$. Note that when $n_0 \ge\lfloor\gamma n_1/(1 - \gamma
) \rfloor+ 1$, $M = \lfloor\gamma n_1/(1 - \gamma) \rfloor+ 1$, and
hence $\max_{1 \le i \le M} m(i) = n_1$. Combining (\ref{equ34}) with the above
fact, we have
\begin{eqnarray}
C_{1,n,\mathrm{SD}}^{(1)} & = & \max_{1 \le n_0 \le n } \max
_{1
\le
i \le M} \biggl\{ \frac{n_0 \alpha_{i+m(i)}^{\prime}}{i} \biggr\} = \alpha.
\nonumber
\end{eqnarray}

Second, we prove that $C_{k, n, \mathrm{SU}}^{(1)} = \alpha$ when $k=1$.
Note that for $i = 1, \ldots, n_0$,
%
\begin{equation}\label{equ35}
\bigl\lfloor\gamma\tilde{m}(i) \bigr\rfloor+ 1 \le\bigl\lfloor\gamma m^*(i)
\bigr\rfloor+ 1 \le i.
\end{equation}
Thus,
%
\begin{equation}\label{equ36}
\frac{n_0 \alpha'_{\tilde{m}(i)}}{i} = \frac{n_0( \lfloor\gamma
\tilde{m}(i) \rfloor+ 1 ) \alpha}{i(n + \lfloor\gamma\tilde
{m}(i) \rfloor+ 1 - \tilde{m}(i))} \le\frac{n_0 \alpha}{n + i -
\tilde{m}(i)} \le\alpha.
\end{equation}
Here, the first inequality follows from (\ref{equ35}) and the second follows
from the fact $\tilde{m}(i) \le i + n_1$. In addition, it is
easy to see that when $i = \lfloor\gamma n \rfloor+ 1$ and
$i + n_1 = n$, we have $m^*(i) = n$ and $n_0 = i$. Thus,
$\tilde{m}(i) = n$ and $\lfloor\gamma\tilde{m}(i) \rfloor+ 1
= i$. By
using the~first equality of (\ref{equ36}), $n_0 \alpha'_{\tilde{m}(i)}/i =
\alpha$. Combining (\ref{equ36})
and the above fact, we have
\begin{eqnarray}
C_{1,n,\mathrm{SU}}^{(1)} & = & \max_{1 \le n_0 \le n } \max
_{1
\le i \le n_0} \biggl\{ \frac{n_0
\alpha_{\tilde{m}(i)}^{\prime}}{i} \biggr\} = \alpha.
\nonumber
\end{eqnarray}\upqed
\end{pf}

\begin{pf*}{Proof of (\ref{equ11})}
As in proving Lemma~\ref{le3.2},
\begin{eqnarray*}
&& \operatorname{Pr}(\widehat R_{n_0} \ge k )
\\
&&\quad = \operatorname{Pr} \Biggl(\bigcup
_{v=k}^{n_0} \{ \widehat P_{(v)} \le
\beta_v \} \Biggr) = \sum_{i=1}^{n_0}
\sum_{r=k}^{n_0} \frac{ \operatorname{Pr}  (\widehat R_{n_0}= r,
\widehat P_{i} \le\beta_r  )}{r}
\\
&&\quad  =  \sum_{i=1}^{n_0} \sum
_{r=k}^{n_0} \frac{ \operatorname{Pr}
(\widehat R_{n_0-1}^{(-i)}= r-1, \widehat P_{i} \le\beta_r  )}{r}
\\
&&\quad  =  \frac{\alpha}{n_0} \sum_{i=1}^{n_0}
\Biggl\{ \sum_{r=k}^{n_0} \operatorname{Pr} \bigl(\widehat R_{n_0-1}^{(-i)} \ge r-1|\widehat P_{i} \le
\beta_{r} \bigr)   -
\sum_{r=k}^{n_{0}-1} \operatorname{Pr}
\bigl(\widehat R_{n_0-1}^{(-i)} \ge r|\widehat P_{i} \le
\beta_{r} \bigr) \Biggr\}
\\
&&\quad \le \frac{\alpha}{n_0} \sum_{i=1}^{n_0}
\Biggl\{ \sum_{r=k}^{n_0} \operatorname{Pr} \bigl(\widehat R_{n_0-1}^{(-i)} \ge r-1|\widehat P_{i} \le\beta
_{r} \bigr)   -
\sum_{r=k}^{n_{0}-1} \operatorname{Pr}
\bigl(\widehat R_{n_0-1}^{(-i)} \ge r|\widehat P_{i} \le
\beta_{r+1} \bigr) \Biggr\}
\\
&&\quad =  \frac{\alpha}{n_0} \sum_{i=1}^{n_0} \operatorname{Pr} \bigl(\widehat R_{n_0-1}^{(-i)} \ge k-1|\widehat P_{i}
\le\beta_{k} \bigr),
\end{eqnarray*}
where the first inequality follows from (A.3) and (A.4) of \citet{r34} and the second follows from Assumption~\ref{ass2b}.
\end{pf*}

\begin{pf*}{Proof of (\ref{equ20})}
Consider a single-step
test based on the $p$-values $\widehat P_1, \ldots,\break  \widehat P_{n_0}$ and
the constant threshold $t$. Let $\widehat R_1$ denote the number of
rejections. Then we have for each $i = 2, \ldots, n_0$,\vspace*{-1pt}
\[
I (\widehat P_{(i)} \le t )  \le I \bigl( \widehat R_1 ( \widehat R_1 - 1 ) \ge i(i-1) \bigr)
\le \frac{1}{i(i - 1)} \sum_{j=1}^{n_0}
\sum_{l (\neq
j)=1}^{n_0} I ( \widehat P_j \le
t, \widehat P_l \le t ),
\]
which proves the desired inequality.
\end{pf*}
\end{appendix}

\section*{Acknowledgements}
Our sincere thanks go to two referees and
the Associate Editor whose very helpful and insightful comments and suggestions
have significantly improved the presentation of the paper.

\begin{supplement}[id=suppA]
\stitle{Supplement to ``Further results on controlling the false discovery~proportion''}
\slink[doi]{10.1214/14-AOS1214SUPP} 
\sdatatype{.pdf}
\sfilename{aos1214\_supp.pdf}
\sdescription{Due to space constraints, we have
relegated to the supplemental article [\citet{r71}] the
remaining figures generated from the simulations in Section~\ref{sec4} and the
findings of additional simulations mentioned in Remark~\ref{re2.1}.}
\end{supplement}



\printaddresses

\end{document}